\font\Bigtenrm=cmr10 scaled\magstep {5.5}
\font\Bigtenbf=cmbx10 scaled\magstep {2}

\vskip 18mm
\centerline{\Bigtenbf G14DIS: Mathematics Dissertation}
\vskip 7mm
\centerline{ Corrected Version $\mit 6$th June $\mit 2000$}
\vskip 2cm
\centerline{\Bigtenrm Algebraic Extensions of Normed Algebras}
\vskip 13cm

\centerline{\sevenrm Disclaimer: This dissertation does not contain plagiarised material; except where otherwise stated all theorems
are the author's.}
\centerline{\sevenrm Acknowledgement: Many thanks to Joel Feinstein for guidance with the literature, useful suggestions and
comments on this work.}
\vskip 2cm
\centerline{by Thomas Dawson}
\vskip 1cm
\centerline{\bf Supervisor: Dr. J. F. Feinstein}
\vfill

\eject

\font\bigtenrm=cmr10 scaled\magstep {2.5}
\font\Bigtenbf=cmbx10 scaled\magstep {3.5}
\def\leaderfill{\leaders\hbox to 1em{\hss.\hss}\hfill}

\noindent {\Bigtenbf Contents}
\vskip 11mm

\noindent
{\bf Chapter 1: Basic Definitions and a Discussion of the Main Problem\hfill 1}\vskip 4pt
\line{1. Introduction\leaderfill 1}
\line{2. Basic Definitions \leaderfill 1}
\line{3. The Main Problem\leaderfill 2}

\vskip 9pt
\noindent
{\bf Chapter 2: Arens-Hoffman Extensions\hfill 4}\vskip 4pt
\line{1. The Main Theorem\leaderfill 4}
\line{2. Maximal Ideal Spaces \leaderfill 8}
\line{3. Discriminants and Tractability\leaderfill 10}

\vskip 9pt
\noindent
{\bf Chapter 3: Examples and an Application of Arens-Hoffman Extensions\hfill 14}\vskip 4pt
\line{1. Uniform Algebras\leaderfill 14}
\line{2. Completions of Normed Algebras and Tractability\leaderfill 15}

\vskip 9pt
\noindent
{\bf Chapter 4: Cole's Construction\hfill 19}\vskip 4pt
\line{1. Algebraic Extension of Uniform Algebras\leaderfill 19}
\line{2. Comparison of Extensions\leaderfill 20}
\line{3. The Cole-Feinstein Conjecture\leaderfill 21}

\vskip 9pt
\noindent
{\bf Chapter 5: Integral Closure of Banach Algebras\hfill 25}\vskip 4pt
\line{1. Introduction\leaderfill 25}
\line{2. Standard Extensions \leaderfill 25}
\line{3. Integral Closure\leaderfill 28}
\line{4. Maximal Ideal Spaces \leaderfill 29}

\vskip 9pt
\noindent
{\bf Chapter 6: Integral Closure of Uniform Algebras\hfill 33}\vskip 4pt
\line{1. Introduction\leaderfill 33}
\line{2. Inverse and Direct Systems \leaderfill 33}
\line{3. Universal Root Algebras \leaderfill 36}

\vskip 9pt
\noindent
{\bf Chapter 7: Conclusion\hfill 37}\vskip 4pt
\line{1. Conclusions\leaderfill 37}
\line{2. Context\leaderfill 37}
\line{3. Further Questions\leaderfill 37}

\vskip 9pt
\noindent
{\bf Appendix 1:} Assumed Gelfand Theory\hfill {\bf 38}

\vskip 9pt
\noindent
{\bf Appendix 2:} Gelfand Theory for Normed Algebras\hfill {\bf 39}

\vskip 9pt
\noindent
{\bf Appendix 3:} Set Theory\hfill {\bf 43}

\vskip 9pt
\noindent
{\bf References \hfill 45}

\vfill\eject


\font\bb=msbm10

\font\eulerfk=eufm10

\font\bbroman=msbm10
\font\smallbbroman=msbm7
\font\tinybbroman=msbm5
\textfont5=\bbroman \scriptfont5=\smallbbroman
 \scriptscriptfont5=\tinybbroman
\def\blackboardroman{\fam5 \bbroman}



 \def\CN{{\blackboardroman C}}	
 \def\NN{{\blackboardroman N}}	
 \def\RN{{\blackboardroman R}}	
 \def\QN{{\blackboardroman Q}}	
 \def\FN{{\blackboardroman F}}	



\def\gamx{\sum\nolimits_{k=0}^pc_kx^k}
\def\delx{\sum\nolimits_{k=0}^qd_kx^k}
\def\betx{\sum\nolimits_{k=0}^{n-1}b_kx^k}

\def\prin{\bigl(\alpha (x)\bigr)}

\def\sqr#1#2{{\vcenter{\hrule height.#2pt
	\hbox{\vrule width.#2pt height#1pt \kern #1pt
		\vrule width.#2pt}
	\hrule height.#2pt}}}
\def\square{\mathchoice\sqr77\sqr77\sqr77\sqr77}

\font\bigtenrm=cmr10 scaled\magstep {2.5}
\font\bigtenbf=cmbx10 scaled\magstep {1.5}

\def\eop{\hfill $\square$\vskip10pt}
\def\thm{\vskip 10pt\noindent THEOREM }
\def\dfn{\vskip 10pt\noindent DEFINITION }
\def\ex{\vskip 10pt\noindent EXAMPLE }
\def\cor{\vskip 10pt\noindent COROLLARY }
\def\lem{\vskip 10pt\noindent LEMMA }
\def\prop{\vskip 10pt\noindent PROPOSITION }
\def\pf{\par\noindent {\it Proof.} }
\def\res{$\ast$}

\def\xb{\bar x}
\def\pri{^\prime}
\def\npp{\nu^{\prime\prime}}
\def\np{\nu^\prime}
\def\Ac{\tilde A}

\def\ker{{\rm ker}\,}	
\def\im{{\rm im}\,}	
\def\OA{\Omega(A)}	
\def\ha{\hat a}		
\def\hA{\hat A}		
\def\eh#1{\varepsilon_#1}	
\def\al{\alpha(x)}
\def\a{\alpha}
\def\ga{\gamma(x)}
\def\Aa{A_{\alpha}}
\def\BA#1{B_#1\colon A_#1}
\def\nls#1{(#1, \Vert\cdot\Vert_{#1})}
\def\vnls#1,#2{(#1_#2, \Vert\cdot\Vert_{#2})}
\def\norm#1{\Vert#1\Vert}
\def\vnorm#1{\left\Vert#1\right\Vert}
\def\vab#1{\left\vert#1\right\vert}
\def\map#1#2{\colon#1\to#2}
\def\mapto#1#2#3#4{\colon#1\to#2;\,#3\mapsto#4}
\def\st{\,\colon\;}
\def\fa#1#2{\forall\,#1\in#2\;}
\def\te#1#2{\exists\,#1\in#2\;}
\def\GU{{\eulerfk U}}
\def\AC{A_\a}


\def\mapright#1{\smash{
	\mathop{\longrightarrow}\limits^{#1}}}

\def\mapup#1{\Big\uparrow
	\rlap{$\vcenter{\hbox{$\scriptstyle#1$}}$}}
\def\mapne#1{\nearrow${$\scriptstyle#1$}$}
\def\mapnw#1{ \nwarrow${$\scriptstyle{#1}$}$  }

\def\maprightt#1#2{\smash{\mathop{\longrightarrow}\limits^{#1}_{#2}}}

\def\dl#1{{\lim_\rightharpoondown}\mathstrut_#1}	
\def\il#1{{\lim_\leftharpoondown}\mathstrut_#1}		


\pageno=1

\font\bigtenrm=cmr10 scaled\magstep {2.5}
\font\bigtenbf=cmbx10 scaled\magstep {1.5}
\noindent {\bigtenbf Chapter 1}
\vskip 7pt
\noindent {\bigtenrm Basic Definitions and a Discussion of the Main Problem}
\vskip 15mm

\noindent {\bf 1. Introduction}\vskip 10pt

\noindent This dissertation is concerned with adjoining roots of polynomials over normed algebras. It was originally 
intended to investigate the Galois theory of topological algebras but this was quickly found to be too ambitious. 
Nevertheless we borrow ideas from field theory to form definitions and try to obtain some analogies.

First we show how to adjoin the root of one monic polynomial to a commutative normed algebra with unity, $A$. This construction is then iterated to adjoin roots of arbitrary sets of monic polynomials, and then to obtain an extension closed under taking of roots of monic polynomials (integrally closed).

There is a similarity with the chain of extensions of $\QN$ (say) to include a root of a single irrational number (say a
 root of $x^2-2\in\QN[x]$), to a field, $\RN$, containing the roots of an uncountable set of polynomials irreducible
  over $\QN$, to an algebraically closed field, $\CN$. For normed algebras this can be achieved entirely by a certain
   type of extension, introduced by Arens and Hoffman in [{\bf 1}]. The simple extension is then iterated using deeper
    set theory to obtain the `standard normed extensions' of chapter 5 which contain roots of any set of monic polynomials
     over $A$. Finally transfinite methods are used on these extensions to obtain an integrally closed normed extension. To
  illustrate the resemblance:
$$\matrix{\hfill\hbox{\sevenrm algebraically closed }&\CN&&C&\hbox{\sevenrm integrally closed}\hfill\cr
	&\mapup\iota&&\mapup\iota&\cr
	&\RN&\hskip 2cm&B&\hbox{\sevenrm a standard extension}\hfill\cr
	&\mapup\iota&&\mapup\iota&\cr
	\hfill\hbox{\sevenrm a simple extension }&\QN(\sqrt2)&&A_\alpha&\hbox{\sevenrm an Arens-Hoffman extension}\hfill\cr
	&\mapup\iota&&\mapup\iota&\cr
	&\QN&&A&}$$

The context of this work in function analysis is described in chapter 7.

A theme is to determine when various properties (principally completeness and semisimplicity) of the original algebra are preserved by the types of extension discussed.

Along the way we are led to consider, in appendix 2, following suggestions in [{\bf 1}], a Gelfand theory for normed rather than Banach algebras.

An application of Arens-Hoffman extensions not directly connected with algebraic extension is included in chapter 3.

Arens-Hoffman extensions are not the only type of extensions
possible. For example in the category of uniform algebras there is a more natural construction available (due to Cole).
We attempt to compare the two methods in the case of a simple extension. Finally in chapter 6 we show that Cole's method can be used to construct integrally closed uniform algebra extensions.

\vskip 20pt
\noindent {\bf 2. Basic Definitions}\vskip10pt

\noindent An algebra is a ring, $A$, which is also a vector space over some field, $\FN$, such that
for all $a,b\in A,\lambda\in\FN$ we have $\lambda(ab)=(\lambda a)b=a(\lambda b)$. A normed algebra is an algebra over $\RN$ or
$\CN$ with a map $\norm\cdot\colon A\to[0,+\infty)$ such that
$(A,\Vert\cdot\Vert)$ is a normed linear space and for all $a,b\in A\;\norm{ab}\le\norm a\norm b$. If the ring $A$ has a unity, $1$, and $\norm 1=1$ then we call $A$ unital, and we call $A$ commutative if it is commutative as a ring. We assume that the reader is familiar with normed spaces and with the basic facts of Gelfand theory (listed in appendix 1).

	Given an algebra, $A$, and an element $a\in A$ there are at least three senses in which $a$ might be called algebraic:
\item{(i)} there is a non-zero $f(x)\in\FN[x]\,:\;f(a)=0$,
\item{(ii)} there is a non-zero $\gamma(x)\in A[x]\,:\; \gamma(a)=0$,
\item{(iii)} there is $n\in\NN$ and $a_0,\ldots,a_{n-1}\in A:\;
a_0+\cdots+a_{n-1}a^{n-1}+a^n=0$.

	Plainly (i) and (iii) are special cases of (ii) when $A$ has a unity. We shall never mean sense (i).
	If (ii) holds then $a$ is {\it algebraic over $A$}. In (iii), $a$ is said to satisfy an {\it equation of
	integral dependence (over $A$)} or to be {\it integral over $A$}. In fact we shall mainly be concerned with integral dependence in this dissertation.

Given $\gamma(x)\in A[x]$ we can ask if there exists $a\in A$ such that $\gamma(a)=0$. If such an $a$ exists it is
called a {\it root} or {\it zero} of $\gamma(x)$. Sometimes a $\ga$ will have no zeros in $A$ and we try to extend
$A$ so as to include one. To be more precise
 we make the following definition:

\dfn 1.1: Let $A$ and $B$ be algebras. $B$ is an {\it extension} of $A$ if there exists a monomorphism $\theta\map A B$.
$B$ is a {\it unital} extension of $A$ if $A$ has a unity, $1$, and $\theta(1)$ is the unity for $B$.
If $A$ and $B$ are normed then $B$ is an {\it isometric} extension of $A$ if there exists an isometric
monomorphism $\theta\map A B$. We shall say that $B$ is a {\it topological} extension if there is a continuous
monomorphism $\theta\map A B$ which is bounded below
(meaning that for some $m>0$ we have $m\norm a\le\norm{\theta(a)}\;(a\in A)$). In each case $\theta$ is referred to as an
{\it embedding}.\par

We shall sometimes abbreviate `$B$ is a unital extension of $A$' to $B:A$.

\noindent The insistence on isometric embedding is explained by Lemma 1.3 and because we only distinguish
normed algebras up to isometric isomorphism.

\dfn 1.2: Let $A$ and $B$ be algebras and $\theta\map A B$ be an embedding. Let $\ga=\gamx\in
A[x]$ and $b\in B\,:\;\theta(\gamma)(b):=\sum\nolimits_{k=0}^p\theta(c_k)b^k=0$. Then $b$ is a {\it root} of $\ga$.

\lem 1.3 (the embedding lemma [well-known]): Let $A, B$ be normed unital algebras and $\theta\map A B$ be a unital
isometric monomorphism. Then there exists a normed algebra, $C$, and a unital isometric isomorphism $\phi\map B C$
such that the inclusion map, $\iota\map A C$ is an isometric embedding and the following diagram is commutative:
$$\matrix{ B&\mapright\phi&C\cr
	\mapup\theta&\mapne\iota&\cr
	A&&\cr}$$
\pf Omitted. (An exercise in set theory.)\eop

\noindent Clearly there are many variants on this.

\vskip 20pt
\noindent{\bf 3. The Main Problem}\vskip 10pt

\noindent Our first problem is, given an algebra, $A$, and a set of polynomials, ${\eulerfk P}\subseteq A[x]$, to find an
extension algebra, $B$, in the same category as $A$ and such that $\forall\,\al\in P\;\exists\,b\in B\st\alpha(b)=0$.
Here we are identifying $\al$ with its image under the map $A[x]\to B[x]$ induced by the implicit embedding $\theta\colon A\to B$. Usually we
shall make the embedding explicit and write $\theta(\a)(x)$ for this.)

Our second problem is to discuss which properties of the original algebra are preserved in $B$. This may cause us to attempt to construct new types of extension algebra. As hinted above we shall only consider unital algebras and then only attempt to adjoin roots of monic polynomials. We now provide some justification for this.

\vskip 20pt
\noindent{\bf Complexification and Unitisation}\vskip 10pt

\noindent If $A$ is a real normed algebra then it is known (see [{\bf 2}] p. 69) that there exists a complex normed algebra,
$A_\CN$, and an isometric monomorphism $\theta\map A {A_\CN}$ (regarding $A_\CN$ as an algebra over $\RN$).
Similarly, if $A$ has no identity we can embed $A$ into a unital normed algebra, $A+\FN$ ([{\bf 2}] p. 15). Therefore
we could take extensions of algebra made complex and unital in which our set of polynomials have roots.

We are of course neglecting questions such as how the final algebra might behave under changing the order of the composite embedding or how the different methods of complexification or unitisation behave with respect to different types of algebraic extension, but these would distract us too much.

The reader may also wonder why we do not assume that $A$ is complete, in view of the availability of a standard completion to embed into.
The main reason is that we would lose the interesting results of chapter 3.

\vskip 20pt
\noindent{\bf Commutativity}\vskip 10pt

\noindent This is not so easily swept aside except to say that the whole problem seems much more natural in
the commutative context. We can, if the coefficients of all polynomials in $P$ commute, consider extensions of the
(maximal or otherwise) commutative subalgebra of $A$ generated by these and the unity, but it seems difficult treat the problem of
extending the whole algebra suitably.

\vskip 20pt
\noindent{\bf Restriction to Monic Polynomials}\vskip 10pt

\noindent First we make the trivial observation that this is really equivalent to the assumption that all polynomials to be
considered have an invertible leading coefficient. The reason why we do not consider adjoining roots of polynomials with
arbitrary leading coefficients is that a general method would tackle the problem of adjoining inverses to $A$. (I.e. to
adjoin an inverse of $a\in A$ we would look to the unital extension generated by $ax-1$.) We refer the reader to [{\bf 3}]
p. 256 for an article on this subject, which is really another problem. (It was proved by \v Silov that if $a\in A$, a
commutative unital Banach algebra, then there exists a unital extension, $B$, containing an inverse for $a$
if and only if $a$ is not a topological zero divisor. So algebraic extension by some polynomials is not possible. A much more
detailed discussion containing some results on this problem can be found in [{\bf 4}])
\vskip 15pt

\hrule\vskip5pt\noindent
Throughout the rest of this dissertation $A$ will denote a complex, commutative, unital normed algebra. All extensions and homomorphisms will be unital unless otherwise stated.
\vskip5pt\hrule\vskip 20pt

\vfill\eject

\noindent {\bigtenbf Chapter 2}
\vskip 7pt
\noindent {\bigtenrm Arens-Hoffman Extensions}
\vskip 15mm

\noindent Let $R$ be a commutative ring with unity and $\al\in R[x]$ be monic.
Let $\prin$ denote the principal ideal
in $R[x]$ generated
by $\al$.
It is well-known that $\al$ has a root in the factor ring
$B:= R[x]\bigm/\bigl(\alpha (x)\bigr)$  and that 
$\nu^\prime\colon R\to B\,;\>a\mapsto\prin+a$ defines an embedding.
The simplicity of this construction is both appealing and practical. 
In the case when $R=A$ is a normed algebra we hope that there exists a norm
on $B$ so that $B:A$ is an isometric extension. 

\par The purpose of this section is relate the work of Arens and Hoffman in
[{\bf 1}] where it is shown that not only can this be done but that the resulting
extensions have good properties. For example if $A$ is complete then so is
$B$.

\par A second concern is whether $B$ is semisimple\footnote{\dag}{See appendix 1.} if $A$ is a semisimple 
Banach algebra. Arens and Hoffman generalise the notion of semisimplicity to normed
algebras and prove their results for these. They make the following
 definition

\dfn 2.1: Let $A$ be a unital normed algebra. $A$ is {\it tractable} if the intersection
of its closed maximal ideals is $\{0\}$.\vskip 10pt

\noindent The motivation for this is explained in appendix 2. Tractability is a stronger condition
than semisimplicity and when $A$ is complete they are equivalent since then it is well-known that maximal ideals
are automatically closed.\vskip 15pt

\hrule\vskip5pt\noindent
Throughout the rest of this chapter $A$ will denote a commutative unital normed
algebra and $\al=a_0+a_1x+\cdots+a_{n-1}x^{n-1}+x^n$ a polynomial over $A$. $B$ will denote
$A[x]\bigm/\bigl(\alpha (x)\bigr)$. The set of all continuous characters (see appendix 2) on $A$ will be denoted $\OA$.
\vskip5pt\hrule\vskip 20pt

\noindent{\bf 1. The Main Theorem}\vskip 10pt

\noindent Arens and Hoffman's solution is to define a norm on the induced polynomial algebra,
$A[x]$, so that $\prin$ is closed. We then have a norm for $B$ in the standard way (see
for example [{\bf 5}] p. 316):
$$\forall\,\beta(x)\in A[x]\qquad \left\Vert\prin+\beta(x)\right\Vert_B:=\inf_{\gamma(x)\in\left(\al\right)}
\norm{\gamma(x)+\beta(x)}.$$
The reader will remember from the proof of this that the function above still defines a
seminorm if $\prin$ is not closed. Most of the proof below establishes that a norm can be found for which this seminorm is faithful.

\lem 2.2\footnote\ddag{With the convention deg(0)=$-\infty$.}: ([{\bf 6}] p. 128)
Let $R$ be a commutative ring with unity
and $$f(x), g(x)\in R[x];\ f(x)=\sum\nolimits_{j=0}^na_jx^j,\ {\rm with}\ g(x)
=\sum\nolimits_{j=0}^mb_jx^j \hbox{ monic } (b_m=1).$$ Then there are unique $q(x), r(x)
\in R[x]$ such that $${\rm deg}\bigl(r(x)\bigr)<{\rm deg}\bigl(g(x)\bigr)\quad
{\rm and}\quad f(x)=q(x)g(x)+r(x)\eqno (*)$$

\noindent {\it Proof.} ([{\bf 6}] p. 128): If ${\rm deg}\bigl(f(x)\bigr)<{\rm deg}\bigl(g(x)\bigr)$ then
we can take $q(x)=0, r(x)=f(x)$ so we may assume $n\ge m$.\par
If $n=0$ then $f(x)=a_0, g(x)=b_0=1\ne 0$. So we have
$$f(x)=a_0=a_0g(x)+0$$ and the result holds. Let $p>0$ and suppose the
result is true for all $n<p$. Put $f_1(x)=f(x)-a_px^{p-m}g(x)$. Then
${\rm deg}\bigl(f_1(x)\bigr)<p$ so by hypothesis there are $q_1(x), r(x)
\in R[x]\quad({\rm deg}\bigl(r(x)\bigr)<{\rm deg}\bigl(g(x)\bigr))\;\colon$
$$f_1(x)=q_1(x)g(x)+r(x)$$
\noindent$\Rightarrow\quad f(x)-a_px^{p-m}g(x)=q_1(x)g(x)+r(x)$\par
\noindent$\Rightarrow\quad f(x)=q_2(x)g(x)+r(x)$ where $q_2(x)= a_px^{p-m}+q_1(x)$.
$r(x)$ has degree less than $m$
so the existence result now follows from the principle of mathematical induction.

\par To show uniqueness, let $f(x),g(x),q(x),r(x)\in R[x]$ satisfy (*) and $q_1(x), r_1(x)\in R[x]$ also having ${\rm deg}\bigl(r_1(x)\bigr)<{\rm deg}\bigl(g(x)\bigr)\quad
{\rm and}\quad f(x)=q_1(x)g(x)+r_1(x)$. Then
$$\bigl(q(x)-q_1(x)\bigr)g(x)=r_1(x)-r(x).$$ So if $q(x)\ne q_1(x)$ then, since
$g(x)$ is monic, the degree of the right hand side is at least $m$,
a contradiction. So $q(x)=q_1(x)$ and $r_1(x)=r(x)$.\eop

\thm 2.3: There exists a norm on $B$ such that $B$ is an isometric extension
of $A$ and such that $B$ is complete if $A$ is.\par

\noindent {\it Proof.} ([{\bf 1}]): Let $$t>0\quad\colon\quad\Vert a_0\Vert+\Vert a_1\Vert t+\cdots+
\Vert a_{n-1}\Vert t^{n-1}\,\le\ t^n$$ and set for $\gamma(x)=\sum\nolimits_{k=0}^pc_kx^k
\in A[x]$,
$$\Vert\gamma(x)\Vert:=\sum_{k=0}^p\Vert c_k\Vert t^k$$

It is easy to check that this defines a norm on $A[x]$:\par
Let $\gamma(x), \delta(x)\in A[x];\ \gamma(x)=\gamx, \delta(x)=\delx$. Let
$r={\rm max}(p, q)$ and $c_k=d_l=0\ \forall k>p, l>q$. Then
$$\eqalign{\Vert\gamma(x)+\delta(x)\Vert=\sum_{k=0}^r\Vert c_k+d_k\Vert t^k\le
\sum_{k=0}^r\bigl(\Vert c_k\Vert+\Vert d_k\Vert\bigr)t^k
&=\sum_{k=0}^r\Vert c_k\Vert t^k+\sum_{k=0}^r\Vert d_k\Vert t^k\cr
&=\Vert\gamma(x)\Vert+\Vert\delta(x)\Vert.\cr}$$
$A[x]$ is an algebra so we need to check that $\Vert\cdot\Vert$ is also submultiplicative:
$$\eqalign{\left\Vert\gamma(x)\delta(x)\right\Vert
&=\left\Vert\sum_{k=0}^{p+q}\left(\sum_{j=0}^kc_jd_{k-j}\right)x^k\right\Vert\cr
&=\sum_{k=0}^{p+q}\left\Vert\sum_{j=0}^kc_jd_{k-j}\right\Vert t^k\cr
&\le\sum_{k=0}^{p+q}\sum_{j=0}^k\Vert c_j\Vert\Vert d_{k-j}\Vert t^jt^{k-j}\cr
&=\bigl\Vert\gamma(x)\bigr\Vert\bigl\Vert\delta(x)\bigr\Vert.\cr}$$
Let $\lambda\in \CN$. Then $$\eqalign{\Vert\lambda\gamma(x)\Vert=
\left\Vert\sum\nolimits_{k=0}^p\lambda c_kx^k\right\Vert
=\sum\nolimits_{k=0}^p\Vert\lambda c_k\Vert t^k
=\sum\nolimits_{k=0}^p\vert\lambda\vert\Vert c_k\Vert t^k
=\vert\lambda\vert\sum\nolimits_{k=0}^p\Vert c_k\Vert t^k
=\vert\lambda\vert\Vert\gamma(x)\Vert}$$
Finally,$$\eqalign{\Vert\gamma(x)\Vert=0\quad\Leftrightarrow\quad\sum_{k=0}^p\Vert c_k\Vert
t^k=0
&\Leftrightarrow\quad\Vert c_k\Vert=0\quad (k=0, \ldots,p)\cr
&\Leftrightarrow\quad c_k=0\quad (k=0, \ldots,p)\cr
&\Leftrightarrow\quad\gamma(x)=0.\cr}$$

\noindent We form the quotient algebra, $B$, as usual
and define $$\eqalign{\left\Vert \prin+\delta(x)\right\Vert 
&=\inf_{\gamma(x)\in A[x]}\Vert\gamma(x)\alpha(x)+\delta(x)\Vert.\cr}$$
With the special choice of $t$, this seminorm is faithful.
Suppose $\prin+\delta(x)\in B$. By the
lemma, there is a unique $q(x)\in A[x]$ and $\beta(x)=
\sum\nolimits_{k=0}^{n-1}b_kx^k\in A[x]\ \colon\ \delta(x)=q(x)\alpha(x)+\beta(x)$.\par
We show that $\forall\,\gamma(x)\in A[x]\ \Vert \alpha(x)\gamma(x)+\beta(x)\Vert
\ge\Vert \beta(x)\Vert$. It then follows that
$$\eqalign{\left\Vert\prin+\delta(x)\right\Vert=
\inf_{\gamma(x)\in A[x]}\Vert\gamma(x)\alpha(x)+\delta(x)\Vert
&=\inf_{\gamma(x)\in A[x]}\Vert(\gamma(x)-q(x))\alpha(x)+\delta(x)\Vert\cr
&=\inf_{\gamma(x)\in A[x]}\Vert\gamma(x)\alpha(x)+\beta(x)\Vert\cr
&\ge \Vert\beta(x)\Vert.\cr}$$ But $\left\Vert\prin+\delta(x)\right\Vert
=\left\Vert\prin+\beta(x)\right\Vert\le\Vert\beta(x)\Vert$ so we shall have
$$\left\Vert\prin+\delta(x)\right\Vert=\Vert\beta(x)\Vert.$$

(In particular, if $\delta(x)\not\in
\prin$ then $\left\Vert\prin+\delta(x)\right\Vert>0$.)
Let $\gamma(x)\in A[x]$. To make the following writing easier we shall set
$a_n=1,b_j=a_j=c_l=0\quad\forall\,l>p,j>n, b_n=0$. Thus

$$\eqalign{\Vert\alpha(x)\gamma(x)+\beta(x)\Vert
&=\left\Vert\betx+\Bigl(\sum\nolimits_{k=0}^na_kx^k\Bigr)\Bigl(\gamx\Bigr)\right\Vert\cr
&=\left\Vert\betx+\sum\nolimits_{j=0}^{n+p}\Bigl(\sum\nolimits_{k=0}^ja_{j-k}c_k
\Bigr)x^j\right\Vert\cr
&=\left\Vert\sum\nolimits_{j=0}^{n+p}\Bigl(b_j+\sum\nolimits_{k=0}^ja_{j-k}c_k
\Bigr)x^j\right\Vert\cr
&=\sum\nolimits_{j=0}^{n+p}\left\Vert b_j+\sum\nolimits_{k=0}^ja_{j-k}c_k
\right\Vert t^j\cr
&=\Vert b_0+a_0c_0\Vert+\Vert b_1+a_1c_0+a_0c_1\Vert+\cdots+
\Vert b_{n-1}+a_{n-1}c_0+\cdots+a_0c_{n-1}\Vert t^{n-1}\cr
&\qquad+\Vert c_0+a_{n-1}c_1+\cdots+a_0c_n\Vert t^n+
\Vert c_1+a_{n-1}c_2+\cdots+a_0c_{n+1}\Vert t^{n+1}\cr
&\qquad+\cdots+\Vert c_p\Vert t^{n+p}\cr
&\ge \bigl(\Vert b_0\Vert-\Vert a_0\Vert \Vert c_0\Vert\bigr)
+\bigl(\Vert b_1\Vert-\Vert a_1\Vert\Vert c_0\Vert-\Vert a_0\Vert\Vert c_1\Vert\bigr)\cr
&\qquad+\cdots+\bigl(\Vert b_{n-1}\Vert-\Vert a_{n-1}\Vert\Vert c_0\Vert-\cdots-
\Vert a_0\Vert\Vert c_{n-1}\Vert\bigr)t^{n-1}\cr
&\qquad+\bigl(\Vert c_0\Vert-\Vert a_{n-1}\Vert\Vert c_1\Vert-\cdots-
\Vert a_0\Vert\Vert c_n\Vert\bigr) t^n\cr
&\qquad+\bigl(\Vert c_1\Vert-\Vert a_{n-1}\Vert\Vert c_2\Vert-\cdots-
\Vert a_0\Vert\Vert c_{n+1}\Vert\bigr)t^{n+1}+\cdots+\Vert c_p\Vert t^{n+p}\cr
&=\Vert b_0\Vert+\Vert b_1\Vert t+\cdots+\Vert b_{n-1}\Vert t^{n-1}\cr
&\qquad+\Vert c_0\Vert\bigl(t^n-\Vert a_{n-1}\Vert t^{n-1}-\cdots-\Vert a_1\Vert t -\Vert a_0\Vert\bigr)\cr
&\qquad+\cdots+t\Vert c_1\Vert\bigl(t^n-\Vert a_{n-1}\Vert t^{n-1}-\cdots\bigr)+\cdots\cr}$$
We see that terms after $\Vert\beta(x)\Vert=\Vert b_0\Vert+\Vert b_1\Vert t+\cdots+\Vert b_{n-1}\Vert t^{n-1}$ are of the form
$$\eqalign{t^j\Vert c_j\Vert\bigl(t^n-\Vert a_{n-1}\Vert t^{n-1}-\cdots-
\Vert a_{n-l}\Vert t^{n-l}\bigr)
&\ge t^j\Vert c_j\Vert\bigl(t^n-\Vert a_{n-1}\Vert t^{n-1}-\cdots-\Vert a_0\Vert\bigr)\cr
&\ge 0\cr}$$ 
(for some $l\in\{0,\ldots,n\}$).
\def\AH{A[x]\bigm/\prin}
So this really does define a norm, $\Vert\cdot\Vert_B$, on $A[x]\bigm/\prin$. 
In particular the natural epimorphism $\nu\colon A[x]\to\AH$ is continuous so
$$\nu^{-1}(\{0\})=\{\delta(x)\in A[x]\ \colon\ \Vert\prin+\delta(x)\Vert=0\}=\prin$$ is closed.
The unital embedding $$\nu^\prime\colon A\to\AH ;a\mapsto\prin+ax^0$$ is an isometry since
the unique polynomial of degree less than $n$ which is equal to $ax^0$ modulo
$\prin$ is $ax^0$ so that $\Vert\prin+ax^0\Vert=\Vert ax^0\Vert_{A[x]}=\Vert a\Vert$.

It remains to check that $B$ is complete if $A$ is. Let $\biggl(\prin+\gamma^{(m)}(x)
\biggr)_{m\in\NN}$ be a Cauchy sequence in $B$. By the lemma, we may assume that
deg$\bigl(\gamma^{(m)}(x)\bigr)<n\ \forall m\in\NN$. So let
$$\gamma^{(m)}(x)=\sum\nolimits_{k=0}^{n-1}c_k^{(m)}x^k\qquad\forall m\in\NN.$$
Then $$\eqalign{\left\Vert\left(\prin+\gamma^{(l)}(x)\right)-\left(\prin+\gamma^{(m)}(x)
\right)\right\Vert
&=\left\Vert\prin+\left(\gamma^{(l)}(x)-\gamma^{(m)}(x)\right)\right\Vert\cr
&=\left\Vert\gamma^{(l)}(x)-\gamma^{(m)}(x)\right\Vert\cr
&=\sum\nolimits_{k=0}^{n-1}\Vert c_k^{(l)}-c_k^{(m)}\Vert t^k\cr
&\to 0\quad {\rm as}\ l,m\to+\infty.\cr}$$
\noindent $\Rightarrow \left(c_k^{(m)}\right)_{m\in\NN}$ is Cauchy ($k=0,\ldots, n-1$)
\par\noindent $\Rightarrow \exists_{k=0}^{n-1}\ c_k\in A\ \colon
\ \lim_{m\to+\infty}c_k^{(m)}=c_k\quad (k=0, \ldots, n-1)$\par So we clearly
have $\lim_{m\to\infty}\gamma^{(m)}(x)=\sum\nolimits_{k=0}^{n-1}c_kx^k=:\gamma(x)$
whence $$\prin+\gamma^{(m)}(x)\to\prin+\gamma(x)\qquad {\rm as}\ n\to+\infty\ {\rm in}\ B.
\eqno\square$$

\noindent Having chosen a specific value of $t$ for the extension we shall refer to it as the parameter of the extension. The theorem shows that a continuum of parameters is possible.
Although we shall not investigate whether or not the Arens-Hoffman norm topology is the only norm topology on $B$ which extends that on $A$ we stop to prove a related remark from [{\bf 7}] to the effect that
there is only one topological solution specified by the family of Arens-Hoffman norms.

\prop 2.4 ([{\bf 7}]): All the norms constructed for $B$ (indexed by $t>0$ where $t^n\ge\Vert a_0\Vert+\Vert a_1\Vert t+\cdots+
\Vert a_{n-1}\Vert t^{n-1}$) are equivalent.
\pf Let $t_j>0\ (j=1,2)$ with $t_j^n\ge\Vert a_0\Vert+\Vert a_1\Vert t_j+\cdots+
\Vert a_{n-1}\Vert t_j^{n-1}\quad (j=1,2)$ and denote the corresponding norms
on $B$ by $\Vert\cdot\Vert_{j,B}\quad (j=1,2)$. Let $t_1\le t_2$ without loss
of generality. It is enough to show that $\exists\,K_1,K_2>0\ \colon
\ \forall\,\beta(x)\in A[x]\ \left({\rm deg}\bigl(\beta(x)\bigr)<n={\rm deg}\bigl(\alpha(x)\bigr)\right)$
$$K_1\left\Vert\prin+\beta(x)\right\Vert_{2,B}\ \le\ \left\Vert\prin+\beta(x)\right\Vert_{1,B}\ \le
\ K_2\left\Vert\prin+\beta(x)\right\Vert_{2,B}$$ Let $\beta(x)=\betx$. Then
$$\left\Vert\prin+\beta(x)\right\Vert_{1,B}=\sum_{k=0}^{n-1}\Vert b_k\Vert t_1^k
\le\sum_{k=0}^{n-1}\Vert b_k\Vert t_2^k=\left\Vert\prin+\beta(x)\right\Vert_{2,B}$$ so
we can take $K_2=1$. Now let $r=t_1/t_2$ and set $K_1={\rm min}(r^0, r,\,\ldots,r^{n-1})$.
Then $K_1>0$ and $$\left\Vert\prin+\beta(x)\right\Vert_{1,B}=\sum_{k=0}^{n-1}
\Vert b_k\Vert t_2^kr^k\ge K_1\sum_{k=0}^{n-1}\Vert b_k\Vert t_2^k=K_1
\left\Vert\prin+\beta(x)\right\Vert_{2,B}.\eqno\square$$

\vskip 15pt
\noindent{\bf 2. Maximal Ideal Spaces}\vskip 10pt

\noindent We now turn to the remaining results in [{\bf 1}]. These concern the spaces of continuous characters
on $A$ and its Arens-Hoffman extension. We recommend that the reader looks at appendix two for the definition
and motivation for studying these spaces. In particular it is explained there why we are studying the space of closed maximal ideals
rather than all the maximal ideals. (The title of the section is therefore a misnomer, but intended to make the reader think
along the right lines with respect to the application of these results.)

The content of the following lemma is not actually needed but it introduces some useful notation.

\lem 2.5 ([{\bf 1}]): Let $\alpha(x)$ be a non-constant monic polynomial over $A$ and let $A[x]$
be given one of the norms above (an `Arens-Hoffman' norm) with parameter $t$.
Then, letting $\Delta_t$ denote the closed disc centre 0 of radius $t$ in $\CN$, we
have (up to a homeomorphism) $\Omega\left(A[x]\right)=\Omega(A)\times\Delta_t$. Also:
$A[x]$ is tractable if and only if $A$ is tractable.

\pf (We fill in some of the detail behind the indentification made
in [{\bf 1}].)
Set $\Upsilon\colon\Omega(A)\times\Delta_t\to\Omega\left(A[x]\right);(h,\lambda)
\mapsto \Upsilon(h,\lambda)$ where, for $\gamma(x)=\gamx\in A[x]$,
$$\Upsilon(h,\lambda)(\gamma(x)):=h(\gamma)(\lambda):=\sum_{k=0}^ph(c_k)\lambda^k.$$
$\Upsilon$ is well-defined: Let $\gamma(x)\in A[x]\quad 
\bigl(\gamma(x)=\gamx\bigr)$. Then it is routine to verify that $\Upsilon(h,\lambda)$ is a homomorphism.
To show it is continuous:
$$\eqalign{\left\vert \Upsilon(h,\lambda)\bigl(\gamma(x)\bigr)\right\vert
&=\left\vert\sum_{k=0}^ph(c_k)\lambda^k\right\vert\cr
&\le\sum_{k=0}^p\vert h(c_k)\vert\vert\lambda\vert^k\cr
&\le\sum_{k=0}^p\Vert c_k\Vert t^k\quad=\Vert\gamma(x)\Vert.\cr}$$
(We used the fact, which follows from Lemma A2.7, that, as for characters on Banach algebras,
the continuous characters on normed algebras have norms at most 1.)

It is also easily checked that $\Upsilon$ is a bijection; for $H\in\Omega(A[x])$ we have
$\Upsilon^{-1}(H)=\left(H\circ\nu^\prime, H(x)\right)$.

\par Finally we show that $\Upsilon$ is a homeomorphism when
$\Omega(A)\times\Delta_t$ is given the product topology.\par
$A$ is unital and so by A2.6, $\Omega(A)$ is compact in the weak
*-topology. $\Delta_t$ is compact by the Heine-Borel theorem and so by
Tychonoff's theorem $\Omega(A)\times\Delta_t$ is compact. Again from a
result in appendix 2, $\Omega\left(A[x]\right)$ is Hausdorff. It is therefore
sufficient to prove that $\Upsilon$ is continuous.
\par Now $\Omega\left(A[x]\right)$ has the weak topology induced by $A[x]\char'136$
(the Gelfand transform: see appendix 2)
so $\Upsilon$ is continuous if and only if for every $\gamma(x)\in A[x]\quad
\gamma(x)\char'136\circ \Upsilon$ is continuous.\par
Let $\gamma(x)=\gamx~\in~A[x]$.
$$\gamma(x)\char'136\circ \Upsilon\colon (h,\lambda)\mapsto\gamma(x)\char'136 \left(\Upsilon(h,\lambda)\right)
=\Upsilon(h,\lambda)(\gamma(x))=\sum_{k=0}^ph(c_k)\lambda^k.$$
Let $(h_\mu, \lambda_\mu)_{\mu\in {\rm M}}$ be a net in $\Omega(A)\times\Delta_t$
with $(h_\mu, \lambda_\mu)\to (h, \lambda)$. So (by definition
of the product topology) $$h_\mu\to h,\qquad\lambda_\mu\to\lambda$$
$\Omega(A)$ has the weak topology induced by $A\char'136$ so we have that
$$\forall\ a\in A\qquad \hat a(h_\mu)=h_\mu(a)\to \hat a(h)=h(a)$$ In particular,
$$h_\mu(c_k)\to h(c_k)\ (k=0, \ldots,p),\qquad \lambda_\mu\to\lambda\qquad\hbox{in }\Delta_t
\subseteq\CN$$
\noindent$\Rightarrow\quad\sum\nolimits_{k=0}^ph_\mu(c_k)\lambda_\mu^k
=(\widehat{\gamma(x)}\circ \Upsilon)(h_\mu, \lambda_\mu)\to
(\widehat{\gamma(x)}\circ \Upsilon)(h, \lambda)$
 (in $\CN$, by the continuity of addition and multiplication there). So
$\Upsilon$ is continuous.\par

 It remains to prove that $A[x]\ \hbox{ is tractable}\ \hbox{if and only if}\  A\
\hbox{ is tractable}.$ Now $A$ and $A[x]$ are commutative, unital normed algebras
so they are tractable if and only if the Gelfand transforms are injective (Theorem A2-9).
First suppose $A[x]$ is tractable. Let $a\in A$ with $\hat a=0$. Set
$\gamma(x)=a$. Then\par\noindent
for each $H\in\Omega(A[x])\ \exists\ h_0\in\Omega(A),\lambda_0\in\CN:
\ H=\Upsilon(h_0,\lambda_0)$\par\noindent
and so $(\gamma(x)\char'136)(H)=(\gamma(x)\char'136)(\Upsilon(h_0,\lambda_0))
=\Upsilon(h_0,\lambda_0)(a)=h_0(a)=\hat a(h_0)=0.$\par\noindent
whence $\gamma(x)\char'136=0$.\par But $A[x]$ is tractable so
$\gamma(x)=ax^0=0$ so $a=0$. Therefore $A$ is tractable. Conversely let $A$ be
tractable. Let $\gamma(x)=\gamx\in A[x]$ with $\gamma(x)\char'136=0$. By
hypothesis we have that for every $H\in\Omega(A[x])\ H(\gamma(x))=0$. So
$\forall\,h\in\Omega(A), \lambda\in\Delta_t\ \Upsilon(h, \lambda)(\gamma(x))
=\sum\nolimits_{k=0}^ph(c_k)\lambda^k=0$. By the uniqueness of the Taylor
coefficients of the analytic function $\sum\nolimits_{k=0}^ph(c_k)z^k$ on
$\Delta_t^\circ$ we have $h(c_k)=0\ k=0, \ldots, p$. But $h$ is arbitrary and
$A$ is tractable so $c_k=0\;(k=0,\ldots, p)$. So $\gamma(x)=0$ and $A[x]$ is
tractable.\eop

\noindent Our next result of Arens-Hoffman, again quoted directly from [{\bf 1}],
is the characterisation of the character space of the extension. We supply a proof for the convenience of the reader. To
ease notation we shall write $\xb$ for the coset $\prin+x$ from now on. 

\thm 2.6: $$\eqalign{\Omega(B)&=\{H\in\Omega(A[x])\,\colon\,H\left(\al\right)=0\}\cr
		&=\{(h, \lambda)\in\Omega(A)\times\Delta_t\;:\;h(\alpha)(\lambda)=\Upsilon(h,
\lambda)(\alpha(x))=0\}.\cr}$$

\pf We can make a similar identification $\amalg:\Omega(B)\to\Omega(A)\times\Delta_t$. Let
$H\in\Omega(B)$. Define $h=H\circ \nu^\prime$ (recall that $\nu^\prime$ was
the natural homomorphism, $\nu\colon A[x]\to B$, restricted to the constants, and is
isometric). So $h$ is a continuous homomorphism $A\to\CN$.
It is non-zero since $h(1)=H\left(\prin+1\right)=H(1_B)\ne 0$.\par
Let $\lambda=H\left(\xb\right)\ \bigl(=(H\circ\nu)(x)\bigr)$. Then $\lambda\in\Delta_t$ since
$$\vert\lambda\vert\le\left\Vert\xb\right\Vert_B\le\Vert x\Vert_{A[x]}=t.$$
So we have $(h, \lambda)\in\Omega(A)\times\Delta_t$ and moreover
$$H\circ\nu=\Upsilon(h, \lambda)\eqno(\dag)$$ because if $\gamma(x)=\gamx\in A[x]$ then
$$\eqalign{H\left(\prin+\gamma(x)\right)&=H\left(\gamma(\xb)\right)
	=H\left(\sum_{k=0}^pc_k\xb^k\right)\cr
	&=H\left(\sum_{k=0}^p \np(c_k)\nu(x)^k \right)\cr
	&=\sum_{k=0}^p (H\circ\np)(c_k)(H\circ\nu)(x)^k \cr
	&=\sum_{k=0}^p h(c_k)\lambda^k=\Upsilon(h,\lambda)(\gamma(x)).\cr}$$

We now check that $\amalg:H\mapsto(H\circ \nu^\prime,(H\circ\nu)(x))$
is injective and that
$$\im\amalg = \{(h, \lambda)\in
\Omega(A)\times\Delta_t\ :\ \Upsilon(h, \lambda)(\alpha(x))=0\}.$$
$\amalg$ is obviously injective because it determines the image of $H$ on
the generators (over $A$), 1 and $\xb$, of $B$.
\par Let $(h, \lambda)\in{\rm im}\,\amalg;\ \exists H\in\Omega(B):
(h, \lambda)=\bigl(H\circ \nu^\prime,(H\circ\nu)(x)\bigr)$. From (\dag), we
have $\Upsilon\bigl(h, \lambda\bigr)\prin=\bigl(H\circ\nu\bigr)\prin
=H\left(\nu\prin\right)=H(0)=0$.\par\noindent
For the reverse inclusion, suppose $\Upsilon(h, \lambda)\prin=0$ for $(h, \lambda)\in
\Omega(A)\times\Delta_t$. Define
$$H:B\to\CN;\ \prin+\gamma(x)\mapsto \Upsilon(h, \lambda)\bigl(\gamma(x)\bigr).$$
This is well-defined: let $\gamma(x), \delta(x)\in A[x]$ with $\prin+\gamma(x)
=\prin+\delta(x)$.\par\noindent
$\Rightarrow\qquad\gamma(x)-\delta(x)\in\prin$\par\noindent
$\Rightarrow\qquad \Upsilon(h, \lambda)(\gamma(x)-\delta(x))=0$\par\noindent
$\Rightarrow\qquad \Upsilon(h, \lambda)(\gamma(x))=\Upsilon(h, \lambda)(\delta(x))$\par
We check that $H\in\Omega(B)$ and that $\amalg(H)=(h, \lambda)$.
Let $\gamx=\gamma(x)\in A[x]$.
We shall write $\gamma(\xb)$ for $\sum\nolimits_{k=0}^pc_k\xb^k$.
Let $\delta(x)\in A[x], \zeta\in\CN$.
$$\eqalign{H(\delta(\xb)+\zeta\gamma(\xb))&=\Upsilon(h, \lambda)\bigl(
\delta(x)+\zeta\gamma(x)\bigr)\cr
&=\Upsilon(h, \lambda)\bigl(\delta(x)\bigr)+\zeta \Upsilon(h, \lambda)\bigl(\gamma(x)\bigr)\cr
&=H(\delta(\xb))+\zeta H(\gamma(\xb)).\cr}$$
$$\eqalign{H\bigl(\delta(\xb)\gamma(\xb)\bigr)&=\Upsilon(h, \lambda)\bigl(
\delta(x)\gamma(x)\bigr)\cr
&=\Upsilon(h, \lambda)\bigl(\delta(x)\bigr)\Upsilon(h, \lambda)\bigl(\gamma(x)\bigr)\cr
&=H(\delta(\xb))H(\gamma(\xb)).\cr}$$
Let $\beta(x)$ be the monic polynomial of degree less than $n$ for
which $\beta(x)-\gamma(x)\in\prin$ $\left(\Leftrightarrow
\beta(\xb)=\gamma(\xb)\right)$. Then $\bigl\Vert\prin+\gamma(x)\bigr\Vert_B=
\Vert\beta(x)\Vert_{A[x]}$ (see Theorem 2.3) so that
$$\vert H(\gamma(\xb))\vert=\vert H(\beta(\xb))\vert=
\vert \Upsilon(h, \lambda)(\beta(x))\vert\le\Vert\beta(x)\Vert_{A[x]}=\Vert\gamma(\xb)\Vert.$$
So $H$ is continuous.
Finally $H(1)=\Upsilon(h, \lambda)(1)=h(1)\ne0$ so $H$ is a character of $B$.
Moreover, $\amalg(H)=\bigl(H\circ \nu^\prime,(H\circ\nu)(x)\bigr)$ and:
$$\forall\,a\in A\quad (H\circ \nu^\prime)(a)=H(a\xb^0)=\Upsilon(h, \lambda)(ax^0)
=h(a)\qquad\Rightarrow\qquad H\circ \nu^\prime=h,$$
$$(H\circ\nu)(x)=H(\xb)=\Upsilon(h, \lambda)(x)=\lambda.$$
Thus $\Upsilon(h, \lambda)=\amalg (H)$.\par
Finally we show that $\amalg$ is a homeomorphism onto a closed subset of
$\Omega(A)\times\Delta_t$. Again it is enough to show $\amalg$ is continuous since
$\Omega(B)$ is compact and $\OA\times\Delta_t$ is Hausdorff by Lemma A2.6. By the definition of the product topology this is
equivalent to showing that the maps
$$\cases{\Omega(B)\to\OA\,;\;H\mapsto H\circ\np,\cr
\Omega(B)\to\Delta_t\,;\;H\mapsto H(\xb)\cr}$$
are continuous. The latter is simply $(\xb)\char'136$ with a restricted codomain and continuous by definition of the topology on $\Omega(B)$. To show that the former is continuous, let
$(H_\lambda)_{\lambda\in\Lambda}$ be a net in $\Omega(B)$ with $\lim_\lambda H_\lambda
=H\in\Omega(B)$. Let $a\in A$. Then $\widehat{\np(a)}\in\hat B$ so
$$\hat a(H_\lambda\circ\np)= (H_\lambda\circ\np)(a)=\widehat{\np(a)}(H_\lambda)\to\widehat{\np(a)}(H)=
(H\circ\np)(a)=\hat a(H\circ\nu).$$ Since $a$ was arbitrary
$H_\lambda\circ\np\to H\circ\np$ in $\OA$.\eop

\vskip5pt\noindent{\bf 3. Discriminants and Tractability}\vskip 10pt

\noindent For the next result we need to define the discriminant of a polynomial over
an arbitrary commutative ring. The definition below is not standard but according to Arens
and Hoffman `agrees with classical usage'. 

\dfn 2.7: Let $R$ be ring and $f(x), g(x)\in R[x]$. The {\it resultant} of $f(x)$ 
and $g(x)$ is, for $f(x)=\sum\nolimits_{k=0}^na_kx^k, g(x)=
\sum\nolimits_{k=0}^mb_kx^k\;(m,n>0)$,
$${\rm res}(f, g):=\left\vert\matrix{
\overbrace{\matrix{a_n &a_{n-1} &\cdots &a_1\cr
			0 &a_n &\cdots &\cdots\cr
			\cdots &\cdots &\cdots &\cdots \cr
			0 &\cdots &\cdots &0\cr}}^{n}
&\overbrace{\matrix{a_0 &0 &\cdots &0\cr
			a_1 &a_0 &\cdots &0\cr
			\cdots &\cdots &\cdots &\cdots \cr
			\cdots &\cdots &a_1 &a_0 \cr}}^{m}	\cr
\matrix{b_m &b_{m-1} &\cdots &\cdots\cr
	0 &b_m &\cdots &\cdots\cr
	\cdots &\cdots &\cdots &\cdots \cr
	0 &0 &\cdots &b_m\cr}
&\matrix{\cdots &\cdots &\cdots &0\cr
	\cdots &\cdots &\cdots &0\cr
	\cdots &\cdots &\cdots &\cdots\cr
	b_{m-1} &\cdots &\cdots &b_0\cr}\cr}\right\vert$$
It can be shown (see [{\bf 6}] p. 325) that if $R=F$, a field,
res$(f, g)=0$ if and only if $a_n=b_m=0\hbox{ or }
f(x)\hbox{ and }g(x)$ have a non-constant common factor in $F[x]$. In particular
we have the familiar criterion from field theory that $f$ has a repeated root 
(in some extension) if and only if res$(f, f^\prime)=0$. We will only apply
these results over the field $\CN$. Since the first $a_n$ can be taken
outside the determinant, if $a_n\ne 0$ one can apply the criterion for
repeated roots to $a_n^{-1}{\rm res}(f, f^\prime)$ in fields: $f$ has a 
repeated root iff $a_n^{-1}{\rm res}(f, f^\prime)=0$. This motivates the definition

\dfn 2.8: Let $R$ be ring with unity and $f(x)\in R[x]$. The {\it discriminant} of $f(x)
=\sum\nolimits_{k=0}^na_kx^k$ is
$${\rm discr}(f):=
\left\vert
\overbrace{\matrix{1 &a_{n-1} &a_{n-2} &\cdots &a_1\cr
			0 &a_n &a_{n-1} &\cdots &\cdots\cr
			\cdots &\cdots &\cdots &\cdots &\cdots \cr
			0 &\cdots &\cdots &\cdots &\cdots\cr
			n &(n-1)a_{n-1} &\cdots &\cdots &a_1\cr
			0 &na_n &(n-1)a_{n-1} &\cdots &\cdots\cr
			\cdots &\cdots &\cdots &\cdots &\cdots \cr
			0 &0 &\cdots &\cdots &na_n\cr}}^{n}
\overbrace{\matrix{a_0 &0 &\cdots &0\cr
			a_1 &a_0 &\cdots &0\cr
			\cdots &\cdots &\cdots &\cdots \cr
			\cdots &\cdots &a_1 &a_0 \cr
			0 &\cdots &\cdots &0\cr
			a_1 &\cdots &\cdots &0\cr
			\cdots &\cdots &\cdots &\cdots\cr
			(n-1)a_{n-1} &\cdots &\cdots &a_1\cr}}^{n-1}
\right\vert$$
This is $a_n^{-1}{\rm res}(f, f^\prime)$ when $R$ is a field.

\ex 2.9: In the commutative ring $R$ we have, for $ax^2+bx+c\in R[x]$,
$${\rm discr}(ax^2+bx+c)=\left\vert\matrix{1 &b &c\cr 2 &b &0\cr
0 &2a &b\cr}\right\vert=b^2-2(b^2-2ac)=4ac-b^2,$$
which is reassuring (though the 1 might not lie in $R$ here).\par

The following theorem is taken directly from [{\bf 1}].

\thm 2.10: Let $A$ be tractable and $d={\rm discr}(\alpha)$ not be a zero divisor or zero.
Then $B$ is tractable.
\pf See [{\bf 1}].\eop

\noindent The following corollary (Arens-Hoffman) is easy to apply and so useful 
for later examples. We reproduce the proof given in [{\bf 1}] to illustrate the use of discriminants.

\cor 2.11: Let $A$ be tractable and $\alpha(x)=x^n+a_0\in A[x]\,(n~\ge~2)$. Then $B$ is tractable
if and only if $a_0$ is not a zero divisor or zero.
\pf $$d={\rm discr}(\alpha)=
\left\vert
\overbrace{\matrix{1 &0 &\cdots &0\cr
			\vdots &\vdots & &\vdots\cr
			0 &\cdots &\cdots &1\cr
			n &0 &\cdots &\cdot \cr
			\vdots &\ddots & &\vdots \cr
			0 &\cdots &\cdots &n\cr}}^{n}
\overbrace{\matrix{a_0 &0 &\cdots &0\cr
			\vdots & & &\vdots\cr
			0 &\cdots &\cdots &a_0\cr
			\cdot & & &0 \cr
			\vdots & & &\vdots\cr
			0 &\cdots &\cdots &0\cr}}^{n-1}
\right\vert$$
We expand this leftwards from the rightmost $a_0$ to obtain
$${\rm discr}(\alpha)=\left((-1)^na_0\right)^{n-1}.
\overbrace{\left\vert\matrix{n &\cdots &0\cr
			 &\ddots &\cr
			 0 &\cdots &n\cr}\right\vert}^{n}=n^na_0^{n-1}.$$
So if $a_0$ is not a zero divisor, $n^na_0^{n-1}$ is not a zero divisor and
by the theorem above, $B$ is tractable.\par
Conversely suppose $d$ is a zero divisor; there exists $a\in A-\{0\}\st aa_0=0$. We
show that $a\xb\in\cap M(B)$, where $M(B)$ denotes the set of closed maximal
ideal of $B$, so that $B$ is not tractable.\par
Let $H\in M(B)$. By Theorem 2.6, $\exists\,(h, \lambda)\in
\Omega(A)\times\Delta_t\st \amalg(H)=\Upsilon(h, \lambda)$ and
$\Upsilon(h, \lambda)\bigl(\alpha(x)\bigr)=0$. Therefore
$\Upsilon(h, \lambda)(x^n+a_0)=\lambda^n+h(a_0)=0$. Also
$$H(a\xb)=\Upsilon(h, \lambda)(ax)=h(a)\lambda.$$Now
$\forall\,h\in\Omega(A)\quad h(aa_0)=h(a)h(a_0)=0$\par$
\Rightarrow\quad\forall\,h\in\Omega(A)\quad h(a)\lambda^n=h(a).-h(a_0)=0.\quad$ So
$$\cases{\lambda=0 &$\Rightarrow H(a\xb)=0$\cr
	\lambda\ne 0&$\Rightarrow H(a\xb)=h(a)\lambda=h(a)\lambda^n.
	\lambda^{-(n-1)}=0$\cr}$$
$H$ was arbitrary so
$$0\ne a\xb\in\bigcap_{H\in\Omega(B)}{\rm ker }H=\cap M(B)\eqno\square$$

\vskip 10pt\noindent
The last main result of the Arens-Hoffman paper is a result saying that although
the Arens-Hoffman extension obtained may not be tractable, it can always be
embedded in a tractable normed algebra in which $\al$ has a root.\par
	Before we can state the theorem we again need to quote some algebraic facts. It is known that if $F$ is a field and $f(x)=a_0+a_1x+\cdots+a_{n-1}
x^{n-1}+x^n\in F[x]$ then the the sums, $q_k=q_k(a_0,\ldots,a_{n-1})$, of the $k$-th powers 
of roots of $f(x)$ ($k=0,\ldots,n-1$) are polynomials in the coefficients 
of $f(x)$ and therefore lie in $F$ (see [{\bf 6}] Theorem 2.20 p. 139). Given a commutative ring, $A$, and $\alpha(x)\in A[x]$ we define $q_k(a_0,\ldots,a_{n-1})$ to be the corresponding polynomial in $a_0,\ldots,a_{n-1}$. Also set $K(B)=\cap M(B)$ (clearly itself a closed ideal of $B$). Then

\lem 2.12: Let the parameter $t$ for the Arens-Hoffman norm on $B$ be large
enough that $t^k\ge \Vert q_k\Vert\quad(k=0,\ldots,n-1)$. Then the $A$ embeds
isometrically in $B/K(B)$, via the natural epimorphism of $B$ onto this
factor algebra.
\pf Let the natural epimorphism be $\nu_B\colon B\to B/K(B)$. It is shown in [{\bf 1}]
that the homomorphism $\npp:=\nu_B\circ\nu^\prime$ is isometric.\eop

\thm 2.13 ([{\bf 1}]): Let $A$ be a tractable normed 
algebra. Then there is a tractable normed algebra extending $A$
isometrically and in which $\alpha(x)$ has a root. If $A$ is complete (so that tractability is equivalent to semismplicity)
then $B$ is a semisimple Banach algebra.

\pf By the preceding lemma, $\npp$ is an isometric monomorphism into the normed
algebra $C:=B/K(B)$. By Theorem 2.3, $B$ is complete when $A$ is so $C$ will be a
Banach algebra when $A$ is. It is obvious that $\alpha(x)$ has a root in $C$; it is $K(B)+\xb$ since
$$\alpha(K(B)+\xb)=K(B)+\alpha(\xb)=K(B)$$
It therefore remains to check that $C$ is always tractable. That a normed
algebra factored by its intersection of closed maximal ideals is tractable
is true in general and proved in Proposition A2.14.\eop

\vfill\eject

\def\da{A_\Delta}	
\def\k{\kappa}
\def\Sb#1{\check S(#1)}
\def\da{A_\Delta}
\def\e{\varepsilon}

\noindent {\bigtenbf Chapter 3}
\vskip 7pt
\noindent {\bigtenrm Examples and an Application of Arens-Hoffman Extensions}
\vskip 15mm

\noindent {\bf 1. Uniform Algebras}\vskip 10pt

\noindent Before going on to look at general questions about algebraic extensions we introduce some examples of normed
algebras which will make useful illustrations. The material in this section is all
standard. An important class is given by the following:

\dfn 3.1: Let $X$ be a non-empty compact Hausdorff space and $A$ be a subalgebra of $C(X)$, the set of all continuous maps $X\to\CN$. $A$ is a {\it uniform algebra} if
\item{(i)} $A$ separates the points of $X$ (meaning that $\fa{\k_1,\k_2}X\ \k_1\ne\k_2\;\Rightarrow\;\te f A\st f(\k_1)\ne f(\k_2)$), 
\item{(ii)} the constant map, $1\in A$,
\item{(iii)} $A$ is uniformly closed (i.e. closed with respect to the supremum norm $\norm f_\infty:=\sup_{\k\in X}\vert f(\k)\vert$).

A more general definition can be given where $X$ is replaced by a locally compact Hausdorff space (see [{\bf 8}] p. 36) but we shall not need this.

\ex 3.2: Let $\Delta=\{z\in\CN\st \vert z\vert\le1\}$. A standard example of a uniform algebra is the {\it disc algebra},
$\da:=\{f\in C(\Delta)\st\hbox{$f$ is analytic on }\Delta^\circ\}$.
(We use $E^\circ$ to denote the interior of a subset, $E$,
of a topological space.)
\vskip10pt

\noindent An extremely useful object associated to any uniform algebra is defined through

\thm 3.3: Let $A$ be a uniform algebra on $X$. Then there exists a unique minimal closed set, $\Sb A\subseteq X$ with the property that every $f\in A$ attains its supremum on $\Sb A$.
\pf We refer the reader to [{\bf 8}] p. 37.\eop

\dfn 3.4: The set $\Sb A$ is called the {\it \v Silov boundary} of $A$.

\ex 3.5: For the disc algebra, $\Sb A\subseteq \partial\Delta$ by the maximum modulus principle. On the other hand,
let $\theta\in(-\pi, \pi]$ and set $f(z)=\left(1\over{2}\right)(1+e^{-i\theta}z)$. Then $f\in A$ but it is easy to see
that $\vab{f(z)}=1=\norm f$ only occurs for $z=e^{i\theta}$. Therefore $\Sb A=\partial\Delta$.
\vskip 10pt
\noindent The next standard result will frequently be referred to.

\thm 3.6: Let $A$ be a uniform algebra on $X$ and $\e\mapto X\OA\k{\eh\k}$ where $\fa\k X\ \eh\k\mapto A\CN f{f(\k)}$. Then $\e$ is a homeomorphism onto a closed subset of $\OA$.
\pf The map is well-defined: let $\k\in X$.
$$\fa{f,g}A,\lambda\in\CN\qquad\eqalign{
&\eh\k(f+\lambda g)=(f+\lambda g)(\k)=f(\k)+\lambda g(\k)=\eh\k(f)+\lambda\eh\k(g)\cr
&\eh\k(fg)=(fg)(\k)=f(\k)g(\k)=\eh\k(f)\eh\k(g)\cr
&\eh\k(1)=1\ne0,\cr}$$
\noindent so $\eh\k\in\OA$.

$\e$ is continuous: Let $(\k_\lambda)_{\lambda\in\Lambda}$ be a net in $X$ with $\k_\lambda\to\k\in X$. Let $f\in A$. Since $f$ is continuous,
$$\hat f( \e_{\k_\lambda} )=\e_{\k_\lambda}(f)=f(\k_\lambda)\to f(\k)=\hat f( \e_\k ).$$
Since the topology on $\OA$ is the one induced by $\hA$ this shows that $\e_{\k_\lambda}\to\e_\k$ as required.

Let $\k,\k\pri\in X$ and suppose $\eh\k=\e_{\k^\prime}$. This means that
$$\fa f A\qquad f(\k)=\eh\k(f)=\e_{\k\pri}(f)=f(\k\pri).$$
Since $A$ separates the points of $X$ this implies that $\k=\k\pri$. So $\e$ is injective. But $\OA$ is Hausdorff (see A1.4) and $X$ is compact so the result follows.\eop
\vskip 10pt\noindent
We therefore have that each $f\in A \hat f$ (which is defined on $\OA\buildrel\e\over\hookleftarrow X$) extends $f$ in the sense that $\hat f\circ\e=f$.

\dfn 3.7:  A uniform algebra on $X$ is {\it natural} (on $X$) if the above map, $\e$, is surjective.

\ex 3.8: The disc algebra, $\da$, is natural. Let $\omega\in\Omega(\da)$. Set $v=\omega(z)$ where $z={\rm id}_\Delta$.
Then $\omega=\eh v$. To prove this we quote without proof the standard result that $\da=P(\Delta)$,
the uniform closure of the algebra ${\rm Pol}_\CN(\Delta)$, the polynomial maps $\Delta\to\CN$ with complex coefficients.
So for all $f\in A\;\exists\,(p_n)\subseteq {\rm Pol}_\CN(\Delta)\st p_n\to f\;(n\to+\infty)$. Now $\omega$ is continuous so $\omega(f)=\lim_{n\to+\infty}\omega(p_n)$.

But for $p(z)=\sum_{k=0}^m a_kz^k\in{\rm Pol}_\CN(\Delta)$ we have
$$\omega(p(z))=\sum_{k=0}^ma_k\omega(z)^k=p(\omega(z)).$$
Thus
$$\omega(f)=\lim_{n\to+\infty}p_n(\omega(z))=f(v)=\eh v(f).$$
\vskip 10pt\noindent
Finally one last standard property of uniform algebras we shall need:

\prop 3.9: Let $A$ be a uniform algebra on $X$. Then $A$ is semisimple.\footnote*{See appendix 1.}
\pf By Theorem 3.6, $X\buildrel\e\over\hookrightarrow \OA$ so by A1.3,
$$\bigcap M(A)\subseteq\bigcap_{\kappa\in X}\ker\eh\k=\{f\in A\st\fa\k X\eh\k(f)=f(\k)=0\}=\{0\}.\eqno\square$$

\vskip 20pt\noindent {\bf 2. Completions of Normed Algebras and Tractability}\vskip 10pt

\noindent The novel condition of tractability featured heavily in the
Arens-Hoffman paper but does not seem to be commonly used in the literature.
One might at first think that this is
because completing\footnote\dag{See appendix 2.} a tractable normed algebra gives a semisimple Banach algebra. But
we shall show (drawing on results in chapter 1) that this is not
always true. This is our first application of Arens-Hoffman extensions to questions other than of algebraic
extension. This work is due to Lindberg.

There are also lots of other questions which spring to mind. Dr. Feinstein
was immediately led to ask whether the operations of forming an Arens-Hoffman
extension of a normed algebra and taking its completion commute. This is stated
in [{\bf 9}] for quadratic Arens-Hoffman extensions; we now prove the general case which has apparently
not been discussed elsewhere. The first part of Lemma 3.10 is a purely algebraic result.

\def\thp{\theta^\prime}

\lem 3.10: Arens-Hoffman extensions satisfy the following universal property. Suppose
$\theta\map {A_2} {B_2}$ is a unital algebra extension and that $A_1$ is another unital algebra. Let $\alpha_1(x)\in
A_1[x]$ be monic and $B_1=A_1[x]/\left(\alpha_1(x)\right)$. Let $\phi_0$ be a homomorphism
$A_1\to A_2$ and $y\in B_2$ be a root of $\alpha_2:=\phi_0(\alpha)(x)$. Then there is a unique homomorphism $\phi\map {B_1}{B_2}$ such that
$$\matrix{B_1&\mapright{\phi}&B_2\cr
	\mapup{\np}&&\mapup{\theta}\cr
	A_1&\mapright{\phi_0}&A_2\cr}$$
is commutative and such that $\phi(\xb)=y$.

If $\BA 1,\BA 2$ are normed and unital with $B_1$ having an Arens-Hoffman norm then $\phi$ is continuous if $\theta,\phi_0$ are continuous.

\pf Let $n={\rm deg}\,\left(\alpha_1(x)\right)$. Recall that every coset
of $B_1$ has a unique representive of degree less than $n$. For
$\beta(x)=\betx\in A_1[x]$ define
$$\eqalign{\phi\left(\beta(\xb)\right)&:=\sum\nolimits_{k=0}^{n-1}(\theta\circ\phi_0)(b_k)y^k\cr
	&=\thp(\beta)(y) \cr}$$
where $\thp:=\theta\circ\phi_0$.
Thus $\phi\circ\np=\theta\circ\phi_0$, and in particular $\phi$ maps the
unity of $B_1$ to the unity of $B_2$.\par
Let $\beta(x)=\betx\in A_1[x],\gamma(x)=\sum\nolimits_{k=0}^{n-1}c_kx^k
\in A_1[x]$, and $\lambda\in\CN$. We have
$$\eqalign{\phi(\beta(\xb)+\lambda\gamma(\xb))
	&=\phi\left(\sum\nolimits_{k=0}^{n-1}\big(c_k+\lambda d_k\big)\xb^k\right)\cr
	&=\sum\nolimits_{k=0}^{n-1}\thp(c_k+\lambda d_k)y^k\cr
	&=\sum\nolimits_{k=0}^{n-1}\left(\thp(c_k)+\lambda\thp(d_k)\right)y^k\cr
	&=\phi(\beta(\xb))+\lambda(\phi(\gamma(\xb))\cr}$$
To show that $\theta$ is multiplicative remember from Lemma 2.1 that
there is a unique
$r(x)\in A_1[x]\,:\ {\rm deg}(r(x))<n \hbox{ and }\exists\,
q(x)\in A_1[x]:$
$$r(x)=\beta(x)\gamma(x)+q(x)\alpha_1(x)\eqno(\dag)$$
Applying the natural epimorphism, $A_1[x]\to B_1$, to (\dag) we obtain $r(\xb)=\beta(\xb)\gamma(\xb)$. Thus $\phi(\beta(\xb)\gamma(\xb))=\phi(r(\xb))=\thp(r)(y)$.

On the other hand $\thp$ induces a homomorphism of the polynomial algebras:
$$\thp_\#\mapto {A_1[x]}{B_2[x]}{\delta(x)}{\thp(\delta)(x)}.$$
Applying $\thp_\#$ to (\dag) gives $\thp(r)(x)=\thp(\beta)(x)\thp(\gamma)(x)+\thp(q)(x)\thp(\a_1)(x)$. Following this by the evaluation-at-$y$ homomorphism gives $\thp(r)(y)=\thp(\beta)(y)\thp(\gamma)(y)
=\phi(\beta(\xb))\phi(\gamma(\xb))$ as required, since we assumed $\theta(\phi_0(\a_1))(y)=0$.

Uniqueness of $\phi$ is obvious since $1, \xb$ generate $B_1$ over $A$.\par
If $A_1, A_2, B_2$ are normed and $\theta, \phi_0$ are continuous (bounded), let $t>0$ be a suitable parameter for the Arens-Hoffman norm on $B_1$. By assumption $\theta$ is also continuous so $\thp$ is continuous. 
For arbitrary $\beta(x)=\betx\in A_1[x]$ we have
$$\eqalign{\Vert\phi(\beta(\xb))\Vert
	&=\left\Vert\sum\nolimits_{k=0}^{n-1}\thp(b_k)y^k\right\Vert\cr
	&\le\sum\nolimits_{k=0}^{n-1}\Vert\thp(b_k)\Vert \norm y^k\cr
	&\le M\sum\nolimits_{k=0}^{n-1}\Vert b_k\Vert t^k\cr
	&=M\Vert\beta(\xb)\Vert,\cr}$$
where $M=\Vert\thp\Vert\hbox{max}_{k=0,\ldots,n-1}\left({\norm y\over t}\right)^k>0$. So $\phi$ is continuous.\eop
 
\eject

\cor 3.11: Let $A_1, A_2$ be commutative unital algebras and
$\phi_0$ be an algebra homomorphism $A_1\to A_2$. Let $\alpha_1(x)
\in A_1[x]$ be monic, $\alpha_2(x)=\phi_0(\alpha_1)(x)$, and $\np_i$ be the canonical embeddings
$A_i\to B_i:=A_i[x]/\left(\alpha_i(x)\right);a\mapsto
\left(\alpha_i(x)\right)+a\quad (i=1,2)$\footnote*{Strictly we should use different symbols for the
inderminates of the two polynomial algebras.}. Then there is a unique
algebra homomorphism $\phi:B_1\to B_2$ with $\phi(\xb)=\xb$ and
such that the following diagram commutes:
$$\matrix{B_1&\mapright{\phi}&B_2\cr
	\mapup{\np_1}&&\mapup{\np_2}\cr
	A_1&\mapright{\phi_0}&A_2\cr}$$
If $A_1, A_2$ are normed (unital), $B_1, B_2$ being given Arens-Hoffman
type norms, and if $\phi_0$ is continuous, then $\phi$ is continuous.\eop

\lem 3.12: Let $A_i\ (i=1,2)$ be commutative unital normed algebras
and $\phi_0$ be an isometric monomorphism $A_1\to A_2$. Let $\alpha_1(x)\in A_1[x]$ be a
monic polynomial, $\alpha_2(x)=\phi_0(\alpha_1)(x)$. Then we can form the
Arens-Hoffman extensions $B_i=A_i[x]/\left(\alpha_i(x)\right) (i=1,2)$
of $A_i (i=1,2)$ using the same norm parameter. Moreover if
$\phi$ is the unique homomorphism extending $\phi_0$ in the above sense
and sending $\xb\mapsto\xb$ then
$\phi$ is an isometry and has dense image if $\phi_0$ does.
\pf Let  $t>0$ satisfy $t^n\ge\sum\nolimits_{k=0}^{n-1}\Vert a_k\Vert t^k$
where $\alpha_1(x)=a_0+\cdots+a_{n-1}x^{n-1}+x^n$. Now
$\alpha_2(x)=\phi_0(a_0)+\cdots+\phi_0(a_{n-1})x^{n-1}+x^n$ and
$\sum\nolimits_{k=0}^{n-1}\Vert \phi_0(a_k)\Vert t^k=
\sum\nolimits_{k=0}^{n-1}\Vert a_k\Vert t^k\le t^n$ so the same parameter for
the norm of $B_2$ can be used. In the same way, $\phi$ must be an isometry:
$\forall\,\beta(x)=\betx\in A_1[x]$ we have
$$\Vert\phi(\beta(\xb))\Vert=
\left\Vert\sum\nolimits_{k=0}^{n-1}\phi_0(b_k)\xb^k\right\Vert
=\sum\nolimits_{k=0}^{n-1}\Vert\phi_0(b_k)\Vert t^k
=\sum\nolimits_{k=0}^{n-1}\Vert b_k\Vert t^k
=\Vert\beta(\xb)\Vert.$$
Now suppose $\im\phi_0$ is dense in $A_2$. Let $\varepsilon>0$, and
$\sum\nolimits_{k=0}^{n-1} c_k^\prime\xb^k\in B_2$. Let $M=
\hbox{max}_{k=0,\ldots,n-1}(t^k)$. By assumption, $\exists\,c_0,\ldots,
c_{n-1}\in A_1\,\colon\ \Vert \phi_0(c_k)-c_k^\prime\Vert<{\varepsilon\over{nM}}
\quad(k=0,\ldots,n-1)$. Thus
$$\eqalign{\left\Vert\phi\Big(\sum\nolimits_{k=0}^{n-1}c_k\xb^k\Big)-
	\sum\nolimits_{k=0}^{n-1} c_k^\prime\xb^k\right\Vert
&=\left\Vert\sum\nolimits_{k=0}^{n-1} (\phi_0(c_k)-c_k^\prime)\xb^k\right\Vert\cr
&=\sum\nolimits_{k=0}^{n-1}\Vert\phi_0(c_k)-c_k^\prime\Vert t^k\cr
&<\sum\nolimits_{k=0}^{n-1}\left({\varepsilon\over{nM}}\right)t^k\ \le\ \varepsilon,\cr}$$
which shows that $\im\phi$ is dense in $B_2$.\eop

We are now ready to prove:	

\thm 3.13: Let $\Aa$ be an Arens-Hoffman extension of $A$ where $\al$ is a monic polynomial over the commutative unital normed algebra $A$ as usual. Then with the same
Arens-Hoffman norm parameter, $t$, we can form the Arens-Hoffman extension of
the completion, $A^\sim$, of $A$. Let this be denoted
$(A^\sim)_\alpha$. Moreover $(A^\sim)_\alpha=(A_\alpha)^\sim$.
\pf Let $(A^\sim, i)$ be a realisation of the completion of $A$. By the lemma
above we can form the Arens-Hoffman extension
$(A^\sim)_\alpha:=A^\sim[x]\bigm/\bigl(i(\alpha)(x)\bigr)$ using the same
value of $t$ for the norm function. Moreover, by Lemma 3.12, $i$ extends to an isometric
monomorphism $\Aa\to(A^\sim)_\alpha$ with dense image. The result follows from
the uniqueness of completions.\eop

\noindent Having estabished Theorem 3.12 it makes sense to consider which
implications hold between the following statements:
\item{(a)} $A$ is tractable
\item{(b)} $(A_\alpha)^\sim$ is semisimple.\par
Answers to these can be found from results in [{\bf 9}] and [{\bf 7}]. We first indicate
why (b) $\Rightarrow$ (a).

\thm 3.14: Let $\alpha(x)$ be a monic polynomial over $A$ and suppose
$(A_\alpha)^\sim$ is semisimple. Then $A$ is tractable.
\pf By Theorem 3.13 we have that $(A^\sim)_\alpha$ is semisimple. It will be
proved in chapter 5 (Corollary 5.8, though it is obvious from Theorem 2.6, or see [{\bf 7}] for a stronger result) that this
implies that $A^\sim$ is semisimple. It now follows from Proposition A2.6 that $A$ is tractable.\eop
\noindent Finally we give an example to show (a) $\not\Rightarrow$ (b).

\thm 3.15 ([{\bf 9}]): There exists a tractable normed algebra, $A$, and a monic polynomial, $\al\in A[x]$ such that $(\Aa)^\sim$ is not semisimple.
\pf The key result used here is Corollary 2.11. The idea of the proof (due to Lindberg who gives a general class of counterexamples in [{\bf 9}]) is to choose $A$ and $\al=x^2-a_0\in A[x]$ so that on completion $A$ remains tractable but $a_0$ becomes a zero divisor so that $(A^\sim)_\a$ is not tractable (which is equivalent to semisimplicity in Banach algebras).

Let $\da$ denote the disc algebra. Set $K=[{1\over 2}, 1], J=K\cup\{0\}$. $\da$ is tractable by Proposition 3.9. Define
$$f\colon J\longrightarrow\CN\,;\,w\longmapsto\cases{
	{1\over w} & $w\ne0$\cr
	0	& $w=0$.\cr} $$
\noindent $f$ is continuous on the compactum, $J$, so by the Stone-Weierstrass theorem, there exists
a sequence, $(f_n)$, of (real) polynomials which converges uniformly to $f$ on $J$.

Consider the subalgebra $A:=\da|_J\,(:=\{g|_J\st g\in\da\})$. Its closure, $\bar A$, in $C(J)$ is a realisation of $\tilde A$.

For $w\in J$ we have $\eh w\in\Omega(\bar A)$. As in Proposition 3.9 we have $\cap M(\bar A)\subseteq\bigcap_{w\in J}\ker\eh w=\{0\}$ so that $\bar A$ is semisimple. The same observation shows that $A$ is tractable.

Now let $f_0$ be the inclusion map $J\to\CN$. We have
$$f_nf_0\longrightarrow\cases{
	1 & uniformly on $K$,\cr
	0 & uniformly on $\{0\}$.\cr} $$
\noindent So $f_nf_0\to e_0\in\bar A\;(n\to+\infty)$ where $e_0$ is the characteristic function of $K$ (in $J$).
We have $e_1:=1-e_0\in\bar A$. From this definition we read off that $e_1f_0=0_J$. Since $e_1, f_0\ne 0$, $f_0$ is a
zero divisor in $\bar A$. Hence by Corollary 2.11, $\bar A[x]/(x^2-f_0)$ is not tractable. This means (since it is complete
by Theorem 2.3) that it is not semisimple.

But by Theorem 3.13 we have $\bar A[x]/(x^2-f_0)=\left(A[x]/(x^2-f_0)\right)^\sim$.\eop

\cor 3.16 ([{\bf 9}]): The completion of a tractable normed algebra need not be semisimple.
\pf In the above example, $A[x]/(x^2-f_0)$ is tractable because $f_0$ is not a zero divisor in $A$. (To see this suppose $f_0g_0=0$ for $g_0\in A$. Let $g_0=g|_J$ where $g\in \da$. We have $(zg)|_J=0$ whence $zg=0$ by the identity principle. So $g$ is zero on $\Delta-\{0\}$ and therefore everywhere (by continuity).)\eop

\vfill\eject

\def\gm#1#2{#1_0^{(\alpha)}+\cdots+#1_{n(\alpha)-1}^{(\alpha)}#2^{n(\alpha)-1}+#2^{n(\alpha)}}
\def\AC{A^\GU}
\def\XC{X^\GU}
\def\p*{\pi^*}
\def\pa{p_\alpha}
\def\ab#1{\vert#1\vert}
\def\da{A_\Delta}

\noindent {\bigtenbf Chapter 4}
\vskip 7pt \noindent {\bigtenrm Cole's Construction}
\vskip 15mm

\noindent {\bf 1. Algebraic Extension of Uniform Algebras}\vskip 10pt

\noindent We describe a construction for adjoining roots of monic polynomials to a uniform algebra.
It was formulated by Cole in [{\bf 10}] for square roots but generalises easily to arbitrary monic polynomials. The
generalisation is very slight and appears to be known to Lindberg (it is implicit in [{\bf 11}]) although the details do
not appear in the literature.
It is interesting for us because it shows we do not have to leave the category of uniform algebras to adjoin the roots
whereas Arens-Hoffman extensions only guarantee the existence of commutative Banach algebras with the desired properties.

Perhaps more immediately striking is that this method allows for roots of arbitrary sets of monic polynomials
to be added all at once.

\vskip 15pt
\hrule\vskip5pt\noindent
Throughout the rest of this chapter $A$ will denote a uniform algebra on the compact Hausdorff space, $X$. $\GU$ will be a
fixed, non-empty set of monic polynomials over $A$. When using $\al$ as an index we will drop the indeterminate.
Let each $\al\in\GU$ be given by $
\al=\gm f x$ for $\left(n(\alpha)\right)_{\alpha\in\GU}\subseteq\NN$ and $f_0^{(\a)},\ldots,f_{n(\a)-1}^{(\a)}\in A$.
\vskip5pt\hrule\vskip 20pt

\thm 4.1: There exists a uniform algebra, $\AC$, and an isometric monomorphism,
$\pi^*\map A \AC$ such that $\fa\al\GU\te f\AC\st \pi^*(\a)(f)=0$.
\pf (Cf. [{\bf 10}].) Let $Y=X\times\CN^\GU$ and give $Y$ the product topology. Let
$$\eqalign{&\pi\map Y X\cr
&p_\a\mapto Y \CN {(\kappa, \lambda)} {\lambda_\a}\cr}$$
be the canonical projections. Let $\pi^*$ be the induced homomorphism
$A\to{C(Y)}\,;\,f\mapsto{f\circ\pi}$. Since $\pi$ is surjective $\p*$ is isometric. Set
$$\eqalign{\XC:
&=\left\{y\in Y\st\fa\a\GU \Bigl(\p*(\a)(\pa)\Bigr)(y)=0\right\}\cr
&=\left\{y\in Y\st\fa\a\GU \Bigl(\p*(f_0^{(\a)})+\cdots+\p*(f_{n(\a)-1}^{(\a)})\pa^{n(\a)-1}+\pa^n\Bigr)(y)=0\right\}\cr
&=\left\{(\kappa, \lambda)\in Y\st\fa\a\GU\
f_0^{(\alpha)}(\kappa)+\cdots+f_{n(\alpha)-1}^{(\alpha)}(\kappa)\lambda_\a^{n(\alpha)-1}+\lambda_\a^{n(\alpha)}=0\right\}. \cr}$$
By the fundamental theorem of algebra (and the axiom of choice) $\XC$ is non-empty; for every $\kappa\in X$ we can choose $\lambda_\a\in\CN\st
f_0^{(\alpha)}(\kappa)+\cdots+f_{n(\alpha)-1}^{(\alpha)}(\kappa)\lambda_\a^{n(\alpha)-1}+\lambda_\a^{n(\alpha)}=0\ (\a\in\GU)$. $\XC$ is Hausdorff because it is a subspace of a Hausdorff space (products of Hausdorff spaces are Hausdorff).

To show $\XC$ is compact we require a lemma. Suppose $c_0,\ldots,c_{m-1}\in\CN, \lambda\in\CN,
\ab\lambda\ge 1$,\par\noindent $c_0+\cdots+c_{m-1}\lambda^{m-1}+\lambda^m=0$.
$\hbox{By assumption,}\qquad -\lambda^m = c_0+\cdots+c_{m-1}\lambda^{m-1}$.\par
$\hbox{Hence}\qquad \ab\lambda^m =\left\vert c_0+\cdots+c_{m-1}\lambda^{m-1}\right\vert$.\par
$\hbox{From which}\qquad \ab\lambda^m \le \ab {c_0}+\cdots+\ab {c_{m-1}} \ab\lambda^{m-1}
	\le \ab {c_0}\ab\lambda^{m-1}+\cdots+\ab {c_{m-1}}
				\ab\lambda^{m-1}\qquad(\hbox{since} \ab\lambda\ge1)$.\par
$\hbox{And so}\qquad \ab\lambda \le \ab {c_0}+\cdots+\ab{c_{m-1}}$.\par
Hence $\fa{(\kappa, \lambda)}\XC \fa\a\GU\ \ab{\lambda_\a}\le 1\hbox{ or }\ab{\lambda_\a}\le 
\left\vert\sum\nolimits_{k=0}^{n(\a)-1} f_k^{(\a)}(\kappa)\right\vert\le
\sum\nolimits_{k=0}^{n(\a)-1}\left\Vert f_k^{(\a)}\right\Vert$. Hence
$$\XC\subseteq X\times\prod_{\a\in\GU}\Delta_{{\rm max}\left(1,\,
\norm{f_0^{(\alpha)}}+\cdots+\norm{f_{n(\alpha)-1}^{(\alpha)}}\right)}$$
which is compact by Tychonoff's theorem. Now by definition, $\XC$ is the intersection of a family of zero sets of continuous functions and so closed. It follows that $\XC$ is compact.

From now on we shall now write $\pi^*, \pa$ for their restrictions to $\XC$.

Now define $\AC$ to be the closed subalgebra of $C(\XC)$ generated by
$\p*(A)\cup\{\pa\st\a\in\GU\}$. So $\p*$ is an embedding into this uniformly closed subalgebra of $C(\XC)$. Clearly $\pa$ is a root of $\al$ for each $\al\in\GU$. So to complete the proof it remains to check that $\AC$ separates the points of $\XC$ but this is obvious because $A$ separates the points of $X$.\eop

\noindent We include one result on the general properties of such extensions. It is explained in [{\bf 10}], at least for square roots, that the operations of passing to the maximal ideal space of a uniform algebra and forming a Cole type of extension commute. We prove a special case here.

\thm 4.2: Let $A$ be natural. Then $\AC$ is natural on $\XC$.

\pf (cf. [{\bf 8}] p. 195) Let the canonical embeddings $X\to\OA$, $\XC\to\Omega(\AC)$ be denoted by $\varepsilon$, $\eta$ respectively. We must show that $\eta$ is surjective, given that $\varepsilon$ is. Let
$\varphi\in\Omega(\AC)$. Consider $\omega:=\varphi\circ\p*\map A \CN$. It is a unital homomorphism since it is a composition of them; in other words $\omega\in\OA$. By hypothesis
$\te \kappa X\st \omega=\eh\kappa$. Thus
$$\fa h A\qquad \varphi(\p*(h))=h(\kappa)$$
Now $\fa\a\GU\ \p*(f_0^{(\a)})+\cdots+\p*(f_{n(\a)-1}^{(\a)})\pa^{n(\a)-1}+\pa^n=0$, so
applying $\varphi$ to this equation gives
$$f_0^{(\alpha)}(\kappa)+\cdots+f_{n(\alpha)-1}^{(\alpha)}(\kappa)\varphi(\pa)^{n(\alpha)-1}+\varphi(\pa)^{n(\alpha)}=0$$
Hence $y:=(\kappa, \lambda)\in\XC$ where $\lambda=\left(\varphi(\pa)\right)_{\a\in\GU}$.

Let $g\in\AC$.There is a sequence of polynomials with complex coefficients in the generating elements, $\p*(A)\cup\{\pa\st\a\in\GU\}$, with $g_n\to g\;(n\to+\infty)$. Explicitly let
$g_n=G_n\left(\p*(h_1),\ldots,\p*(h_{l(n)}),p_{\a_1},\ldots,p_{\a_{l(n)}}\right)$ for a non-decreasing sequence
$\left(l(n)\right)_{n\in\NN}\subseteq\NN$ and $(h_n)\subseteq A, (\a_n)\subseteq\GU$. Then
$$\eqalign{\varphi(g)&=\varphi\left(\lim_{n\to+\infty}g_n\right)\cr
&=\lim_{n\to+\infty}\varphi(g_n)\cr
&=\lim_{n\to+\infty}G_n\left(h_1(\kappa),\ldots,h_{l(n)}(\kappa),\varphi(p_{\a_1}),\ldots,\varphi(p_{\a_{l(n)}})\right)\cr
&=\lim_{n\to+\infty}G_n\left(\p*(h_1)(y),\ldots,\p*(h_{l(n)})(y),p_{\a_1}(y),\ldots,p_{\a_{l(n)}}(y)\right)\cr
&=\lim_{n\to+\infty}g_n(y)=g(y).\cr}$$
So $\varphi=\eta_y$.\eop

\vskip 20pt\noindent {\bf 2. Comparison of Extensions}\vskip 10pt

\noindent We would now like to compare the two methods of extending a uniform algebra so as to inlude a root of the polynomial $\al\in A[x]$. To do this systematically we modify the appropriate definitions from field theory.

\dfn 4.3: Let $A_1, A_2, B_1, B_2$ be commutative algebras with unities and $\theta_k
\map {A_k} {B_k}$ be unital monomorphisms ($k=1,2$). Let $\phi_0\map {A_1} {A_2}$ be a unital
isomorphism. Then $B_1\colon A_2$ and $B_2\colon A_2$ are {\it isomorphic extensions} if there exists a unital isomorphism $\phi\map {B_1} {B_2}$ making the following diagram commute:
 $$\matrix{B_1&\mapright{\phi}&B_2\cr
	\mapup{\theta_1}&&\mapup{\theta_2}\cr
	A_1&\mapright{\phi_0}&A_2\cr}$$
If $A_1=A_2=A$ and $\phi_0={\rm id}_A$ then we shall say that $B_1, B_2$ are {\it isomorphic over $A$} and write $B_1\cong_A B_2$.\par
If $A_1, A_2, B_1, B_2$ are normed unital algebras and $\theta_1, \theta_2$ are continuous and $\phi_0, \phi$ are homeomorphisms, $\BA 1, \BA 2$ will be called {\it topologically isomorphic extensions}. If all maps are isometries then $\BA 1, \BA 2$ will be said to be {\it isometrically isomorphic extensions}.

\def\Ap{A^\prime}
\def\Xp{X^\prime}
\noindent All the remaining results in this chapter are apparently new.

\ex 4.4: Let $\da$ be the disc algebra and $\al=x^2-z$ where $z={\rm id}_\Delta$. Let $\Ap$ denote the extension generated by Cole's method; $\Ap$ is the closed subalgebra of $C(\Xp)$ generated by $\p*(\da)\cup\{\pa\}$ in the notation from above where $\Xp=
\{(u,v)\in\Delta\times\CN\st v^2=u\}=\{(v^2, v)\st v\in\Delta\}$.

Let $\Aa$ denote the corresponding Arens-Hoffman extension. Lemma 3.10 shows that $$\phi\mapto \Aa \Ap
{f+g\xb} {\p*(f)+\p*(g)\pa}$$ is a homomorphism whose image is clearly a dense subspace of $\Ap$. It will follow very
easily from a result in the next section that $\phi$ is a topological isomorphism $\Aa\to\Ap$ over $\da$ (the hard part of this is to show that $\phi$ is surjective).

We could therefore renorm $\Aa$ so as to make $\phi$ an isometry. But the resulting norm will not be an Arens-Hoffman
norm for any choice of norm parameter $t>0$.

To see this, suppose $t>0$ is such that $\norm{f+g\xb}=\norm f + \norm g t$ defines a norm on $\Aa$ so that $\phi$ is an isometry. We must have
$$t=\norm\xb =\norm {\phi(\xb)}=\norm \pa=\sup_{(v^2, v)\colon v\in\Delta}\left\vert\pa(v^2, v)\right\vert
=\sup_{v\in\Delta}\ab v =1.$$
But now let $f=({1\over 2})(z+1), g=({1\over 2})(z-1)$. We must have
$$\eqalign{\norm {f+g\xb}&=\left\Vert \p*(f)+\p*(g)\pa\right\Vert\cr
	&=\sup_{v\in\Delta}\left\vert \bigl(\p*(f)+\p*(g)\pa\bigr)(v^2, v) \right\vert\cr
	&=\sup_{v\in\Delta}\left\vert f(v^2)+vg(v^2) \right\vert\cr
	&=\left\vert f(e^{2i\theta})+e^{i\theta}g(e^{2i\theta})\right\vert\cr}$$
for some $\theta\in(-\pi, \pi]$ by the maximum modulus principle. Thus
$$\norm {f+g\xb}\le\left\vert f(e^{2i\theta})\right\vert+\left\vert g(e^{2i\theta})\right\vert
=(1/2)\left(\bigl\vert e^{2i\theta}+1\bigr\vert +\bigl\vert e^{2i\theta}-1\bigr\vert \right)\qquad (\dag)$$
$$\hbox{Now }\cases{\vert e^{2i\theta}+1\vert<2 &unless $\theta=0,\pi$, and\cr
			\vert e^{2i\theta}-1\vert <2 &unless $\theta=\pm{\pi\over 2}$,\cr}$$
so (\dag) shows that $2=\norm f+\norm g=\norm {f+g\xb}<(2+2)/2=2$, a contradiction.\eop

\def\gm#1#2{#1_0^{(\alpha)}+\cdots+#1_{n(\alpha)-1}^{(\alpha)}#2^{n(\alpha)-1}+#2^{n(\alpha)}}
\def\AC{A^\GU}
\def\XC{X^\GU}
\def\p*{\pi^*}
\def\pa{p_\alpha}
\def\ab#1{\vert#1\vert}
\def\da{A_\Delta}
\def\aA{A^\a}
\def\r*{\rho^*}
\def\hB{\hat B}
\def\chB{(\hB)^-}
\def\Gt{\tilde \Gamma}
\def\Sb{\check S(A)}

\vskip 15pt\noindent {\bf 3. The Cole-Feinstein Conjecture}\vskip 10pt

\hrule\vskip5pt\noindent
In this section $A$ will be a uniform algebra on the compactum $X$ (`$(A, X)$ is a uniform algebra'). We shall write $(A^\a, X^\a)$ for the simple Cole extension $(A^{\{\a\}}, X^{\{\a\}})$. The corresponding Arens-Hoffman extension, $A[x]/\prin$, will be denoted $\Aa$ or $B$ as usual.
\vskip5pt\hrule\vskip 20pt

\noindent It is easy to give examples where $\aA\ne\Aa$. For instance, $\bar A$ of Theorem 3.15 is obviously a uniform
algebra on $J$ but there was a monic polynomial over $\bar A$ for which the Arens-Hoffman extension was not semisimple
so it could not be isomorphic to a Cole type of extension because those are semisimple (Proposition 3.9).

More generally, by Corollary 2.11, $A[x]/(x^2-f_0)$ cannot be isomorphic to a uniform algebra if $f_0\in A$ is a zero divisor or zero.

However Theorem 2.13 states that factoring the Arens-Hoffman extension by its radical, $K(B)$, leads to a semisimple
isometric extension of $A$. (Remember that uniform algebras are complete so $B$ is always complete and so
semisimplicity is equivalent to tractability then.) The proper question to ask is therefore whether or not we always
have $B/K(B)=:C=\aA$. B. J. Cole and J. F. Feinstein conjectured that this is true but we shall show that there are
counterexamples.

\vskip10pt\noindent
Recall that there is a continuous injection $\varepsilon\map X \OA$ and that (modulo $\amalg$)
$$\eqalign{\Omega(B)&=\{(h,\lambda)\in\OA\times\Delta_t \st \Upsilon(h,\lambda)(\al)=0\}\cr
			&=\{(h,\lambda)\in\OA\times\Delta_t \st \sum\nolimits_{k=0}^{n-1}h(a_k)\lambda^k=0\},\cr}$$
where $t>0$ is the Arens-Hoffman norm parameter. We exploit the strong resemblance of this with
$$X^\a=\{(\kappa, \lambda)\in X\times\CN\st a_0(\kappa)+\cdots+a_{n-1}(\kappa)\lambda^{n-1}+\lambda^n=0\}$$
Now $t$ only had to satisfy $t^n\ge\sum\nolimits_{k=0}^{n-1}\norm {a_k}t^k$ so we can assume
$t\ge\max(1,\norm{a_0}+\cdots+\norm{a_{n-1}})$. Therefore by the lemma in the proof of Theorem 4.1 we have that
for every ${(\kappa, \lambda)}\in{X^\a}\ab\lambda\le t$. We summarise:

\lem 4.5: $\rho\mapto{X^\a}{\Omega(B)}{(\kappa, \lambda)}{\amalg^{-1}(\eh\kappa,\lambda)}$ is a homeomorphism onto a closed subspace of $\Omega(B)$.

\pf Since $\amalg$ is a homeomorphism it is enough to check that the map
$(\kappa, \lambda)\mapsto(\eh\kappa,\lambda)$ is a homeomorphism onto a closed subspace of $\OA\times\Delta_t$.
(The map is well defined by the comments above.) The codomain is Hausdorff and $X^\a$ is compact so we only need to
check injectivity and continuity. Since $\varepsilon$ is continuous and injective the map composed with each of the canonical projections onto $\OA, \Delta_t$ has these properties and the result follows.\eop

\noindent $\rho$ induces a map $\r*\mapto{C(\Omega(B))}{C(X^\a)}{h}{h\circ\rho}$. This map is clearly a (unital)
homomorphism. We first show that it restricts to a map $\hB\to\aA$ where $\hB\,\left(\subseteq C(\Omega(B))\right)$
is the image of the Gelfand transform, $\Gamma\map{B}{\hB}$. By the fundamental isomorphism theorem, $\Gamma$ induces
an isomorphism $\Gt\map{B/K(B)}{\hB}$ of norm $\norm\Gamma\le 1$.

\lem 4.6: $\r*(\hB)\subseteq\aA$.

\pf A typical element of $B$ is $b=\beta(\xb)$ where $\beta(x)=\betx\in A[x]$. We must show that $\r*(\hat b)=\hat b\circ\rho\in\aA$.\par
Let $(\kappa, \lambda)\in X^\a$. Please refer to Theorem 2.6 for the notation used.
$$\eqalign{\left(\widehat{\beta(\xb)}\circ\rho\right)(\kappa, \lambda)
	&=\widehat{\beta(\xb)}\left(\amalg^{-1}(\eh \kappa,\lambda)\right)\cr
	&=\left(\amalg^{-1}(\eh \kappa,\lambda)\right)\left(\beta(\xb)\right)\cr
	&=\left(\amalg^{-1}(\eh \kappa,\lambda)\circ\nu\right)\left(\beta(x)\right)\cr
	&=\Upsilon(\eh \kappa,\lambda)\left(\beta(x)\right)\qquad\hbox{(by Theorem 2.6)}\cr
	&=\sum\nolimits_{k=0}^{n-1}\eh\kappa(b_k)\lambda^k\cr
	&=\sum\nolimits_{k=0}^{n-1}b_k(\kappa)\lambda^k\cr
	&=\sum\nolimits_{k=0}^{n-1}\p*(b_k)(\kappa,\lambda)\pa(\kappa,\lambda)^k\cr
	&=\left(\sum\nolimits_{k=0}^{n-1}\p*(b_k)\pa^k\right)(\kappa,\lambda)\cr}$$
Thus $\r*(\hat b)=\sum\nolimits_{k=0}^{n-1}\p*(b_k)\pa^k\in\aA$.\eop

\noindent We can therefore consider $\r*$ to be a contractive homomorphism $\hB\to\aA$. Moreover since $\aA$ is generated by $\im\p*\cup\{\pa\}$, $\r*$ has dense image. The last line of the above proof allows us to summarise with a commutative diagram:
$$\matrix{B/K(B)&\mapright\Gt&\hB&\mapright{\r*}&\aA\cr
	&\mapnw{\npp}&&\mapne{\p*}&\cr
	&&A&&\cr}$$
There is one immediate positive result:

\prop 4.7: Let $A$ be natural on $X$. Then $\aA=\chB$ (isometric isomorphism).

\pf In this case $\varepsilon$ is a homeomorphism so $\rho$ is surjective. Hence for $b\in B$,
$\norm{\r*(\hat b)}=\norm{\hat b\circ\rho}=\norm{\hat b}$ and $\r*$ is an isometry. Now the unique continuous linear extension of $\r*$ to $\r*\map{\chB}{\aA}$ is therefore also isometric. But $\chB, \aA$ are complete and $\overline{\im\r*}=\aA$ so $\r*$ is surjective (considered as a map $\chB\to\aA$).\eop

\noindent It is therefore desirable to know when $\hB$ is closed in $C(\Omega(B))$ (i.e. complete in the sup-norm).
For natural uniform algebras this has been worked out by Heuer and Lindberg in [{\bf 12}].  (First recall

\dfn 4.8: $a\in A$ is a {\it topological divisor of zero} if $\exists\,(u_n)\subseteq A\st u_na\to0\;(n\to+\infty)$ with
$\fa n\NN\norm{u_n}=1$)

\thm 4.9 ([{\bf 12}] p. 338): Let $A$ be a natural uniform algebra and $B=\Aa$ an Arens-Hoffman extension.
Let $d:=\hbox{discr}(\al)$ not be a zero divisor.
Then $\hB$ is sup-norm complete if and only if $d$ is not a topological divisor of zero.
\vskip 10pt\noindent
We therefore have:

\cor 4.10: Let $A$ be natural and $d$ not be a zero divisor.
Then $\aA=\hB\iff d$ is not a topological divisor of zero.\eop

\noindent To give examples of this situation there is a further equivalent condition available:

\thm 4.11 ([{\bf 13}] p. 137): Let $A$ be a uniform algebra and $d\in A$. Then $d$ is a topological zero divisor
if and only if $d$ has a zero on $\check S(A)$.
\vskip10pt\noindent
Before going on, note that the condition on $d$ in Corollary 4.10 is quite strong; if $d$ is not a
zero divisor then by Theorem 2.10, $B$ is already semisimple. Therefore we can prove the following

\thm 4.12: Let $A$ be a natural uniform algebra. Then $B$ and $\aA$ are topologically isomorphic extensions of $A$ if $d$ is not a topological zero divisor.
\pf Suppose $d$ is not a topological divisor of zero. By Theorem 2.10, $K(B)=\{0\}$ so that $B/K(B)=B$
(to isometric isomorphism). By Proposition 4.7 and Theorem 4.9, $\aA=\chB$ and $\chB=\hB$. By Banach's isomorphism theorem, $\Gt$ is a topological isomorphism $B/K(B)\to\hB$.\eop 

\noindent Thus for example, for the disc algebra, $\da$, and $\alpha(x)=x^2-z$ we have discr$(\al)=4z$ (see example 2.9)
which does not vanish on $\check S(\da)=\partial\Delta$ (example 3.5) so that Cole's extension is topologically isomorphic
to the Arens-Hoffman extension.

\vskip 10pt
\noindent We have not yet found a counterexample to the Cole-Feinstein conjecture although we do now know that if $T$ is a topological isomorphism such that the following commutes
$$\matrix{B/K(B)&\mapright T &\aA\cr
	\mapup\npp&\mapne{\pi^*}&\cr
	A&&\cr}\eqno(\dag)$$
and $T(K(B)+\xb)=\pa$ then $T=\r*\circ\Gt$. This forces $\r*$ to be both injective and surjective. $\r*$ is always a
contraction so if $\hB$ is not uniformly closed (i.e. closed in the sup norm) then (when $A$ is natural) we would have a
continuous algebraic isomorphism $\hB\to\chB$ which seems strange but not impossible.
However for our situation this cannot happen
and I am very grateful to Dr. Feinstein for
pointing out the following result to me which is an elementary exercise in uniform algebras:

\lem 4.13 (Dr. Feinstein): Suppose that $\vnls B,1$ is a commutative semisimple complex unital Banach algebra and that
$(B_2, \norm\cdot_\infty)$ is a uniform algebra. Suppose that $T\map{B_1}{B_2}$ is an algebra isomorphism. Then $B_1$ is
topologically isomorphic to its the Gelfand
transform. In particular $\widehat{B_1}$ is complete in the sup norm.

\def\snorm#1{\left\Vert#1\right\Vert_\infty}
\pf Since $B_2$ is a uniform algebra we have for each $f\in B_2\ \snorm f=\lim_{n\to+\infty} \snorm{f^n}^{1\over n}$. Now $T$
induces a second norm, $\vab\cdot$ on $B_1$ defined by
$$\eqalign{\fa b{B_1}\qquad \vab b &:=\snorm{T(b)}=\lim_{n\to+\infty}\snorm{T(b)^n}^{1\over n}\cr
	&=\lim_{n\to+\infty}\snorm{T(b^n)}^{1\over n}\cr
	&=\lim_{n\to+\infty}\vab{b^n}^{1\over n}.\qquad\qquad (\dag)\cr }$$
\noindent In this way, $T^{-1}$ is an isometric isomorphism $(B_2, \norm\cdot_\infty)\to (B_1,\vab\cdot)$. Since
$B_2$ is complete, $(B_1,\vab\cdot)$ is complete. Hence by the (commutative) uniqueness of norm theorem, $\vab\cdot$ and $\norm\cdot_1$ are
equivalent norms.

Now $\widehat{B_1}$ has the sup norm inherited from $C(\Omega(B_1))$ and by Corollary A1.8 for $b\in B_1$ we have
$$\vnorm{\hat b}=\lim_{n\to+\infty}\vnorm{b^n}^{1\over n}_1.$$

\noindent Since $\vab\cdot$ and $\norm\cdot_1$ are equivalent norms it follows easily from (\dag) that the Gelfand transform is
continuous and bounded below as a map from $(B_1, \vab\cdot)$ and so also from $\vnls B,1$.\eop

\noindent
We now come to our main result allowing for counterexamples to the Cole-Feinstein conjecture to be found very easily.

\thm 4.14: Let $A$ be a natural uniform algebra and $d={\rm discr}(\al)$ be a topological divisor of zero in $A$, but not zero
or a zero divisor. Then there does not exist an (algebraic) isomorphism $B/K(B)\to \aA$.

\pf Suppose $A$ and $d$ satisfy the hypotheses above. Since $d$ is not zero or a zero divisor, $B$ is semisimple by
Theorem 2.10. So $K(B)=\{0\}$ and $B$ is isometrically isomorphic to $B/K(B)$ where the latter space has the quotient norm.
Suppose for a contradiction that there is an isomorphism $T\map B{\aA}$. By Lemma 4.13, $\hat B$ must be complete
and this contradicts Theorem 4.9.\eop

\noindent So for example letting $A=\da$, the disc algebra, and $f_0=z-1, \al=x^2-f_0$ we have that $d=4(z-1)$ is not a zero
divisor or zero but is a topological zero divisor (by Theorem 4.11). So by the last theorem, $B/K(B)$ and $\aA$ are essentially
different here. [By Proposition 4.7 there is however an embedding $B/K(B)\hookrightarrow\aA$.]

We close this chapter by remarking that a condition on $\al$ commonly occurring in the literature
(`separability') on the type of polynomial considered excludes
examples like this one.

\vfill\eject

\def\GA{{\cal U}}
\def\GM{{\eulerfk M}}
\def\GL{{\eulerfk L}}
\def\ao{\alpha_0}
\def\b{\beta}
\def\g{\gamma}
\def\PS{{\cal P}(S)}
\def\CQ{{\cal Q}}
\def\Map{{\rm Map}}
\def\dom{{\rm dom}\,}
\def\pa{p_\alpha}

\noindent {\bigtenbf Chapter 5}
\vskip 7pt
\noindent {\bigtenrm Integral Closure of Banach Algebras}
\vskip 15mm

\noindent {\bf 1. Introduction}\vskip 10pt

\noindent In field theory it is well-known that algebraic closures exist. These are field extensions which contain roots of all polynomials over themselves. In other words there is no need to find further extensions in order to adjoin roots. In the same way one can obtain extensions which are closed under taking roots of monic polynomials. More precisely:

\dfn 5.1: Let $A$ be a commutative algebra with identity. $A$ is {\it integrally closed} if for every monic $\al\in A[x]\;\te\xi A\st \a (\xi)=0$. If $B\colon A$ is an algebra extension and $B$ is integrally closed then $B$ is an {\it integral closure} of $A$.\vskip 10pt

\noindent Note that we do not ask (here) that every monic polynomial `splits' completely into linear factors in an 
integral closure; for results on this problem see [{\bf 7}].
Nor do we ask that integral closure be in any sense minimal here.

The main theorem in this section is that every Banach algebra has an integral closure (in the same category). At the moment we can only adjoin roots of finitely many polynomials at a time so our first task is to extend this method. Lindberg has made a study of this in [{\bf 11}] so we are reporting his results.
\vskip 15mm

\noindent {\bf 2. Standard Extensions}\vskip 10pt

\noindent Let $\GU\subseteq A[x]$ be any non-empty set of monic polynomials where $A$ is a commutative algebra with unity. We follow Lindberg by obtaining an algebra extension in which every polynomial in $\GU$ has a root and then consider the case when $A$ is normed.

\dfn 5.2 ([{\bf 11}]): Let $B\colon A$ be an algebra extension and $\GU\subseteq A[x]$ a set of monic
polynomials. $B\colon A$ is called a {\it standard extension} if there exists a well-ordering, $\le$, on $\GU$ (with $\ao=\inf\GU$) and intermediate subalgebras, $(B_\a)_{\a\in\GU}$ such that
\item{(i)} $B=\bigcup_{\a\in\GU} B_\a$
\item{(ii)} $\fa{\a,\beta}\GU\ \a\le\b\,\Rightarrow\,B_\b\colon B_\a$ (with respect to inclusion)
\item{(iii)} $\fa{\a}\GU B_\a\cong_{A_\a}\,A_\a[x]/\prin$ where
$$A_\a:=\cases{\bigcup_{\b<\a} B_\b& $\ao<\a$\cr
	A & $\ao=\a$.\cr}$$

\noindent (Warning: we are using $\Aa$ to mean something different to in previous chapters.)
Recall that the notation in (iii) means that there is an isomorphism $\phi_\a\map{B_\a}{A_\a[x]/\prin}$
satisfying $\phi_\a(a)=\np_\a(a)\ (a\in A_\a)$ where $\np_\a$ is the usual embedding $\Aa\to\Aa[x]/\prin$.
It is clear that Arens-Hoffman extensions are standard extensions and that the embedding lemma can be used to construct
standard extensions from any finite sequence of Arens-Hoffman extensions. We should also warn the reader that the definition
in [{\bf 11}] is slightly more general in that it allows the polynomial, $\al$, generating each intermediate subalgebra,
$B_\alpha$, to have coefficients in $\Aa$ rather than $A$.

\thm 5.3 ([{\bf 11}]): Let $A$ be a commutative algebra with unity and $\emptyset\ne\GU\subseteq A[x]$ a set of monic
polynomials. Then there exists a standard extension, $B_\GU\colon A$, in which every $\al\in\GU$ has a root.
\pf (We follow the outline in [{\bf 11}] and supply extra details for the convenience of the reader.)
By the well-ordering principle there exists a well-ordering, $\le$, on $\GU$. Let $\ao=\inf\GU$. We shall need to refer to some results from set theory in the course of this proof. These are given in appendix 3 to avoid cluttering this proof, but the reader may wish to glance at it now to check notations and conventions.

First we require a set, $S$, with $\#S>\#A$ and $A\subset S$. (For example ${\cal P}(A)\cup A$ could be used for $S$.) Let $\CQ$ be the set of all maps from initial segments of $\GU$ to $\PS$:
$$\CQ:=\Map(\GU,\PS)\cup\bigcup_{\a\in\GU}\Map\left([\ao, \a), \PS\right).$$
($\GU$ is an initial segment by convention.) We now let
$$\GM:=\left\{\,f\in\CQ\st\eqalign{
\hbox{(i) }&\fa\a{\dom f}\ f(\a)\colon A\hbox{ w.r.t. inclusion}\cr
\hbox{(ii) }&\fa{\a,\b}{\dom f}\ \a\le\b\,\Rightarrow\,f(\b)\colon f(\a)\hbox{ w.r.t. inclusion}
\cr
\hbox{(iii) }&\fa\b{\dom f}\ f(\b)\cong_{\tilde f(\b)}\tilde f(\b)[x]/(\beta(x))\hbox{ where}\cr
\tilde f(\b)&:=\bigcup_{\a<\b}f(\a)\hbox{ if }\ao<\b\hbox{ and }\tilde f(\b):=A\hbox{ if }\ao=\b\cr}\ 
\right\}$$
Note that (i) implies that $f(\a)$ carries an algebraic structure. Really what we have written is not adequate, but it is consistent with the universal abuse of language when speaking of `the algebra $A$'. An algebra is strictly a 4-tuple, $(A, S_+, S_\cdot, S_{\FN\times A})$, where $S_+, S_\cdot\subseteq (A\times A)\times A, S_{\FN\times A}\subseteq(\FN\times A)\times A$ define addition, multiplication and the action of the field on $A$. (If $A$ is normed there is also a fifth component, $S_{\Vert\cdot\Vert}\subseteq A\times\RN_+$ determining the norm map.) So it would be more correct to consider a suitable subset of maps, $f$, into
$$S\times{\cal P}(S\times S\times S)\times{\cal P}(S\times S\times S)\times{\cal P}(\FN\times S\times S),$$ but this level of formality is not appropriate here.

$\GM\ne\emptyset$ for the trivial map, $\emptyset$, the single element of $\Map([\ao, \ao), \PS)$, belongs to $\GM$. (The reader can be assured that the `inductive step' in the proof below will make it clear that there are non-trivial elements in $\GM$.) For $f,g\in\GM$ define
$$f\preceq g\quad:\Leftrightarrow\quad \dom f\subseteq\dom g\,\wedge\, g|_{\dom f}=f.$$
This clearly defines a partial ordering on $\GM$. Suppose $\emptyset\ne\GL\subseteq\GM$ is totally ordered.
Now $\bigcup_{f\in\GL}\dom f=:D$ is an initial segment of $\GU$ for if $D\subset\GU$ we have $D=[\ao, \a)$ where
$\a=\inf(\GU-D)$. In either case we can define $g\map D \PS$ by $\a\mapsto f(\a)$ where $f\in\GL, \a\in\dom f$. This
is well-defined by the assumption on $\GL$ and easily checked to be a member of $\GM$. Since $g$ is an upper bound
for $\GL$, $(\GM, \preceq)$ is inductively ordered. By Zorn's lemma there is a maximal element, $f\in\GM$.

Suppose $\dom f\subset\GU$. We shall show that this leads to a contradiction and then use $f$ to define $B_\GU$.

Let $\a=\inf(\GU-\dom f)$ so that $\dom f=[\ao, \a)$. Set
$$\Aa:=\cases{\bigcup_{\b<\a} f(\b)&$\ao<\a$\cr
	A &$\ao=\a$.\cr}$$
Thus $\Aa$ is an algebra extension of $A$ when it is given the obvious structure (which is well-defined by the conditions in $\GM$).

We claim that $\#\Aa=\#A$. Now since $\GU\subseteq A[x]$ and $\# A[x]=\# A$ (Lemma A3.4) it is enough (Lemma A3.3) to show that
$\forall\,\b<\ao\;\#f(\b)=\#A$. This we do by transfinite induction.

Let $J=\{\b<\a\st \#f(\b)=\#A\}$ and suppose $[\ao, \b)\subseteq J$.

Now $f(\b)\cong\tilde f(\b)[x]/(\beta(x))$ so $\#f(\b)=\#\tilde f(\b)[x]/(\beta(x))$. Recall from chapter 2 that since $\tilde f(\b)$ is
a commutative ring with unity this factor algebra is identifiable setwise with
$\mathop{\times}\limits_{k=0}^{{\rm deg}(\beta(x))-1}\tilde f(\b)$. This last set has cardinality $\#\tilde f(\b)$ by Lemma A3.2
(all our algebras will be infinite). But by definition of $\tilde f(\b)$, the inductive hypothesis, and Lemma A3.3 again, $\#\tilde f(\b)=\#A$. Therefore $\b\in J$ and the claim follows from the principle of transfinite induction.

Now form $B_\a^\prime:=\Aa[x]/\prin$. By the claim and by the same reasoning as for the factor algebra considered just above, $\#B_\a^\prime=\#A$. Let $\np$ be the usual embedding of $\Aa$ into this algebra. We shall need a variant of the embedding lemma to complete the diagram below:
$$\matrix{B_\a^\prime&\cdots\rightarrow&B_\a\subseteq S\cr
	\mapup\np&\mapne\subseteq&\cr
	\Aa&&\cr
	\mapup\subseteq&&\cr
	A&&\cr}$$
By assumption, $\#S>\#A=\#\Aa$ so $\#(S-\Aa)>\#A=\#B_\a^\prime\ge\#(B_\a^\prime-\im\np)$. Let $i$ be an injection
$B_\a^\prime-\im\np\to S-\Aa$. Define the subset $f_1(\a):=\im i\,\dot\cup\Aa$ and let it have the algebraic structure determined by
the induced bijection $B_\a^\prime\to f_1(\a)$. Thus $f_1(\a)$ is a commutative algebra unitally extending $\Aa$.

\vskip 5mm

$$\matrix{&&&&i&&&&\cr
	&&B_\a^\prime-\im\np&\qquad&\longrightarrow&\qquad&\im i&\qquad&\cr
	B_\a^\prime\vphantom{\Biggl|}&&&&&&&&S\cr
	&&\im\np&&\longleftarrow&&\Aa&&\cr
	&&&&\np&&&&\cr}$$

\vskip 5mm

Define $f_1(\b)=f(\b)$ for $\b<\a$. If $\GU=[\ao, \a]$ then $\dom f_1$ is an initial segment of $\GU$. In the case
$\GU\supset[\ao, \a]$ we have $\dom f_1=[\ao, \inf(\GU-[\ao, \a]))$, also an initial segment of $\GU$. It is routine to check that
$f_1$ satisfies (i),(ii),(iii) so that $f_1\in\GM$. However we have $f_1\succ f$, a contradiction.

Therefore $\dom f_0=\GU$ and we define
$$B_\GU:=\bigcup_{\a\in\GU}f(\a)$$
In this way $B_\GU$ is a standard extension of $A$ relative to $(B_\a)_{\a\in\GU}$ where $\fa \a \GU B_\a:=f(\a)$. We have that $\fa \al\GU \al$ has a root in $B_\a\subseteq B_\GU$ by condition (iii).\eop

\noindent So we have solved the algebraic problem. We next consider the normability of such extensions.

\thm 5.4 ([{\bf 11}]): Let $\nls A$ be a commutative unital normed algebra and $B$ a standard extension
of $A$ with respect to $(B_\a)_{\a\in\GU}$ where $\GU\ne\emptyset$ is a well-ordered set of monic polynomials
over $A$. Let the isomorphisms $B_\a\cong_{A_\a}A_\a[x]/\prin$ be $\phi_\a\map{B_\a}{A_\a[x]/\prin}\;(\a\in\GU)$. Then there exists a
norm, $\norm\cdot_B$, on $B$ and values of the Arens-Hoffman norm parameters for the $A_\a[x]/\prin\;(\a\in\GU)$ where the
norm on $\vnls A,\alpha$ is the restriction of $\norm\cdot_B$ to $\Aa$, such that $\phi_\a$ is an isometry for all $\a\in\GU$.

\pf This is an application of the transfinite recursion theorem, which we use informally. 

By assumption the following commutes for all $\a\in\GU$:
$$\matrix{B_\a&\maprightt{\phi_\a}{\cong}&A_\a[x]/\prin\cr
	\mapup\subseteq&\mapne{\np_\a}&\cr
	\Aa&&\cr}$$
where $\np_\a$ is our usual embedding of an algebra into its Arens-Hoffman extension; it will be an isometry whatever value is chosen for the norm parameter. 

Let $\g\in\GU$ and suppose that for all $\a<\g$, $B_\a$ has a norm $\norm\cdot_\a$ and $t_\a>0\;(\a<\g)$ are such that
\item{(i)} $\forall\,\a<\b<\g\ \norm\cdot_\b\hbox{ extends }\norm\cdot_\a$
\item{(ii)} $\forall\,\b<\g\ \phi_\b$ is an isometry when $A_\b$ has norm the restriction of $\norm\cdot_\b$ and $A_\b[x]/\left(\b(x)\right)$ is given an Arens-Hoffman norm with parameter $t_\b$.

Now $A_\g=\bigcup_{\b<\g}B_\b$ so by assumption (i), $\left(\vnls B,\b\right)_{\b<\g}$ induce a well-defined norm on $A_\g$.
Choosing $t_\g>0$ according to the condition in Theorem 2.3 makes $A_\g[x]/\left(\g(x)\right)$ into an Arens-Hoffman extension when given an Arens-Hoffman norm. We can now define a norm, $\norm\cdot_\g$, on $B_\g$ by declaring $\phi_\g$ to be an isometry. By the comments in the first paragraph the system $\left(\vnls B,\b\right)_{\b\le\g}$ satisfies (i), (ii) above but with `$<\g$' replaced with `$\le\g$'. The result now follows.\footnote*{The reader may be more familiar with this type of informal inductive definition when $\GU$ is replaced by an initial segment of some ordinal number. However there is no essential difference in view of result A3.16 of appendix three.}\eop

\vskip 15mm

\noindent {\bf 3. Integral Closure}\vskip 10pt

\noindent We are now in a position to prove the main theorem of this chapter. It is a simplification of a result stated
in [{\bf 11}] and the proof is too.
Similar methods are used in [{\bf 14}] but this has not been our main source.

\thm 5.5: Let $\nls A$ be a commutative unital Banach algebra. Then there is a
Banach algebra extension, $\nls C$, which is integrally closed.
\pf Let $\omega_1$ denote the first uncountable ordinal, and 0 the first
ordinal. Let, for $\gamma\in [0, \omega_1]$, $p(\gamma)$ be the proposition
`There exist Banach algebras, $\vnls C, \alpha_{\alpha<\gamma}$ such that:
\item {(i)}$\qquad \vnls C,0=\nls A$
\item {(ii)}$\qquad \forall\,\alpha\le\beta<\gamma\qquad
\vnls C,\beta\hbox{ is a unital Banach algebra extension of }
\vnls C,\alpha\hbox{ with respect to }\subseteq$
\item{(iii)}$\qquad \forall\,\alpha<\beta<\gamma\quad\forall\,
\hbox{monic }\mu(x)\in C_\alpha[x]\ \exists\,\xi\in C_\beta\,:\ \mu(\xi)=0$.'
\par

Suppose $\gamma<\omega_1$ and that $p(\gamma)$ holds. We must show that we
can choose a Banach algebra, $\vnls C,\gamma$, so that the family
$\vnls C, \alpha_{\alpha\le\gamma}$ satisfies (i)-(iii) above. This is
obvious if $\gamma=0$ so let $0<\gamma$.\par

	Set $A_\gamma=\bigcup_{\alpha<\gamma}C_\alpha$ and give $A_\gamma$ the
[well-defined] algebraic operations and norm,
$\Vert\cdot\Vert_{A_\gamma}$, induced by the subalgebras
$\vnls C, \alpha_{\alpha<\gamma}$.\par

By Theorem 5.4 we can form a standard normed extension, $\nls {B_\gamma}$
of $\nls {A_\gamma}$ in which every monic polynomial over $A_\gamma$ has
a root. Recall that in this construction the unital embedding was given by
inclusion. Therefore every monic polynomial over $C_\alpha\ (\alpha<\gamma)$
has a root in $B_\gamma$.\par

$B_\gamma$ may not be complete so we form its completion, $B_\gamma^\sim$.
By the embedding lemma, there exists a Banach algebra, $\vnls C,\gamma$, containing
$A_\gamma$ as a subalgebra, extending its norm and such that there is an
isometric isomorphism $\phi\colon B_\gamma^\sim\to C_\gamma$ with
$\phi\circ j=\iota$ where $j$ is the embedding of $B_\gamma$ in  $B_\gamma^\sim$
and $\iota$ is the inclusion map $B_\gamma\to C_\gamma$.\par

Let $\alpha<\gamma$ and $\mu(x)\in C_\alpha[x]$ be monic. So $\mu(x)$ is a
monic polynomial over $A_\gamma$ and so there is $\xi\in B_\gamma\,:\,
\mu(\xi)=0$. But $B_\gamma\subseteq C_\gamma$ so $p(\gamma+1)$ is satisfied by 
this family of Banach algebras.\par

Therefore $p(\omega_1)$ holds; let $\vnls C,\alpha_{\alpha<\omega_1}$ be a
family of Banach algebras satisfying $p(\omega_1)$.\par

Set $C=\bigcup_{\alpha<\omega_1}C_\alpha$ and give $C$ the usual induced
algebraic structure and norm. It remains to check that $C$ is complete and
integrally closed.\par

Let $(c_n)\subseteq C$ be a Cauchy sequence. So $\exists_{n\in\NN}\,
\alpha_n<\omega_1\,:\,\forall\,n\in\NN\ c_n\in C_{\alpha_n}$.
As in the proof of Theorem A3.20, $\exists\,\alpha<\omega_1\colon\ 
\forall\,n\in\NN\ \alpha_n<\alpha$. So $(c_n)\subseteq C_\alpha$. But
$C_\alpha$ is complete so $\exists\,c\in C_\alpha\subseteq C\colon\,
c_n\to c\quad(n\to+\infty)$. Thus $C$ is complete.\par

Finally let $\mu(x)$ be a monic polynomial over $C$. It has finitely many
coefficients so $\exists\,\alpha<\omega_1\colon\ \mu(x)\in\ C_\alpha[x]$. But
then $\mu(x)$ has a root in $C_{\alpha+1}\subseteq C$. So $C$ is integrally
closed.\eop
\vskip 15mm

\noindent {\bf 4. Maximal Ideal Spaces}\vskip 10pt

\noindent This chapter could really have been called `standard extensions' and though their main application here was to form an integral closure of an arbitrary Banach algebra they generalise Arens-Hoffman extensions. One can ask if the results of chapter two also generalise. The first thing to note is that the method of obtaining a standard normed extension from a set of monic polynomials, $\GU\subseteq A[x]$, will not generally give a Banach algebra, even if $A$ is Banach. The details of exactly when what has been called $B_\GU$ is Banach (given $A$ is complete) have been worked out in [{\bf 11}].

We turn our attention instead to determining the space of continuous characters (so again our section title is a deliberate
misnomer) of a standard normed extension and are
able to generalise Theorem 2.6 directly. This refines a remark in [{\bf 11}] (for our definitions) that the map
induced by taking restrictions $\Omega(B_\GU)\to\Omega(A)$ is surjective. (The reasoning is different: we
obtain $\Omega(B_\GU)$ explicitly whereas Lindberg makes use of the fact (which there is not space to discuss
here) that standard extensions are examples of `integral extensions'.)

\vskip 15pt
\hrule\vskip5pt\noindent
For rest of this chapter $\GU$ will be a non-empty set of monic polynomials over $A$:
$$\fa\a\GU\qquad\al=a_0^{(\a)}+\cdots+a_{n(\a)-1}^{(\a)}x^{n(\a)-1}+x^{n(\a)}.$$ $\le$ will be a well-ordering on $\GU$ with least element $\ao$.
\vskip5pt\hrule\vskip 20pt

\def\OB{\Omega(B)}

\thm 5.6: Let $B:=B_\GU$ be the generalised Arens-Hoffman extension of the normed algebra $A$ as in Theorem 5.4. Then (up to homeomorphism)
$$\Omega(B)=\left\{(h, \lambda)\in\OA\times\CN^\GU\st\quad\fa\a\GU\quad h(\a)(\lambda_\a)=0\right\}$$

\pf Let $P$ denote the set on the right-hand side of this equation.

Let $t_\a>0$ be the parameter for the Arens-Hoffman norm on $\Aa[x]/\prin\;(\a\in\GU)$. We see from the
inductive construction that it is possible to choose $t_\a\ge 1,\left\Vert a_0^{(\a)}\right\Vert+\cdots+\left\Vert
a_{n(\a)-1}^{(\a)}\right\Vert\ \fa\a\GU$. Let all notation be as in Theorem 5.4, so that $\phi_\a\map{B_\a}{\Aa[x]/\prin}$
is an isometric isomorphism $(\a\in\GU)$. Define
$$\amalg\mapto\OB P H {(h, \lambda)}$$
where $h=H|_A$ and $\fa\a\GU \lambda_\a=H\left( \phi_\a^{-1}(\xb) \right)$. For convenience we shall set $\xi_\a:=\phi_\a^{-1}(\xb)\;(\forall\,\a\in\GU)$. 
This is well-defined since $H|_A$ is a continuous unital homomorphism $A\to\CN$ and for $\a\in\GU$ we have
$$\eqalign{h(\a)(\lambda_\a)&=h(a_0^{(\a)})+\cdots+h(a_{n(\a)-1}^{(\a)})\lambda_\a^{n(\a)-1}+\lambda_\a^{n(\a)}\cr
&=H \left(a_0^{(\a)}+\cdots+a_{n(\a)-1}^{(\a)}\xi_\a^{n(\a)-1}+\xi_\a^{n(\a)}\right)\cr
&=(H\circ\phi_\a^{-1})\left(a_0^{(\a)}+\cdots+a_{n(\a)-1}^{(\a)}\xb^{n(\a)-1}+\xb^{n(\a)}\right)\cr
&=(H\circ\phi_\a^{-1})(0)=0.\cr}$$
An easy induction shows that $B$ is generated (over $A$) by $\{1\}\cup\{\xi_\a\st\a\in\GU\}$.\footnote*{It is not possible that the $\xi_\a$ will not all be distinct so we can use set notation here.} Therefore $\amalg$ is injective.

$\amalg$ is surjective: Let $(h, \lambda)\in P$. We use transfinite recursion (informally) to construct a suitable lift of $h$ to $B$.

Let $\b\in\GU$ and suppose that $(H_\a)_{\a<\b}$ is a family of continuous characters ($\forall\,\a<\b\;H_\a\in\Omega(B_\a)$):
\item{(i)} $\forall\,\a<\b\quad H_\a|_A=h$,
\item{(ii)} $\forall\,\a<\b\quad H_\a(\xi_\a)=\lambda_\a$,
\item{(iii)} $\forall\,\a<\g<\b\quad H_\a=H_\gamma|_{B_\a}$.

Now define $h_\b\mapto{A_\b}\CN b{H_\a(b)}$ where $b\in B_\a$ (or if $\b=\ao$ define $h_\b=h$). This is well-defined by
(iii) and plainly $h_\b\in\Omega(A_\b)$.

By assumption $h( a_0^{(\b)} )+\cdots+h( a_{n(\b)-1}^{(\b)} )\lambda_\b^{n(\b)-1}+\lambda_\b^{n(\b)}=0$. This implies
(as in the proof of Theorem 4.1) that
$$\vert\lambda_\b\vert\le1\hbox{ or }\vert\lambda_\b\vert\le\vab{h(a_0^{(\b)})}+\cdots+\vab{h(a_{n(\b)-1}^{(\b)})}
\le \vnorm{a_0^{(\b)}}+\cdots+\vnorm{a_{n(\b)-1}^{(\b)}}$$
(since $\norm h\le1$ by A2.7). So $\vert\lambda_\b\vert\le t_\b$.

Therefore $(h_\b, \lambda_\b)\in\Omega(A_\b)\times\Delta_{t_\b}$ and by Theorem 2.6, there exists $\omega\in\Omega\left(A_\b[x]/(\beta(x))\right)\st (h_\b, \lambda_\b)=(\omega\circ\np_\b, \omega(\xb))$.

Define $H_\b=\omega\circ\phi_\b$. Thus $H_\b\in\Omega(B_\b)$. It is simple to check that $H_\b$ is consistent with (i), (ii), (iii):
\item{(i)} Let $a\in A$. $H_\b(a)=\omega((\beta(x))+a)=h_\b(a)=h(a)$ so $H_\b|_A=h$
\item{(ii)} $H_\b(\xi_\b)=\omega(\xb)=\lambda_\b$
\item{(iii)} Suppose $\a<\b$. Let $b\in B_\a\subseteq A_\b$. Then $H_\b(b)=\omega(\phi_\b(b))=\omega(\np_\b(b))=h_\b(b)=H_\a(b)$.

By the transfinite recursion theorem\footnote\dag{In this case the objects are sufficiently simple for a formal application of that theorem. In the notation of A3.17 we can take $W=\GU$ and\par $X=\bigcup_{\hbox{\sevenrm initial segments}, J,\hbox{\sevenrm of } \GU}\left\{(H_\a)_{\a\in J}\in\mathop{\times}\limits_{\a\in J}\Omega(B_\a)\st \hbox{statements similar to (i), (ii), (iii) hold}\right\}$. } there exist continuous characters $H_\a\in\Omega(B_\a)\;(\a\in\GU)$ such that
\item{(i)} $\forall\,\a\in\GU\qquad H_\a|_A=h$,
\item{(ii)} $\forall\,\a\in\GU\qquad H_\a(\xi_\a)=\lambda_\a$,
\item{(iii)} $\fa{\a,\b}\GU\qquad \a\le\b\ \Rightarrow H_\a=H_\b|_{B_\a}$.

Therefore we can define $H\mapto{B=\bigcup_{\a\in\GU}B_\a}\CN b{H_\a(b)}$ where $\a\in\GU, b\in B_\a$. Thus $H$ is a well-defined unital homomorphism $B\to\CN$ and bounded since $\vab{ H(b)} =\vab{ H_\a(b) }\le\norm b$. So $H\in\OB$. Moreover we have $\amalg(H)=(h, \lambda)$.

$\OB$ is compact by A2.6 and $P$ is clearly Hausdorff so it remains to check that $\amalg$ is continuous;
the proof of this is similar to the corresponding part of Theorem 2.6.\eop

\noindent Our first application of this is to prove the result mentioned at the start of this section.

\cor 5.7: Let $B:A$ be a standard normed extension relative to $(B_\a)_{\a\in\GU\subseteq A[x]}$. Then the following are surjections:
$$\eqalign{&\OB\to\OA\,;\;{H}\mapsto{H|_A}\cr
&\OB\to{\Omega(B_\a)}\,;\;H\mapsto{H|_{B_\a}}\qquad(\a\in\GU)\cr}$$
\pf This is clear for the first map. For the second, let $H_\a\in\Omega(B_\a)$.
$$\hbox{Let}\qquad\lambda_\b=\cases{H_\b(\xi_\b) & if $\b\le\a$\cr
	\lambda\in\CN\st H_\a(\b)(\lambda_\b)=0 & (for some specific choices) if $\b>\a$.\cr}$$
Thus $H:=\amalg^{-1}(H|_A, \lambda)$ and $H_\a$ agree on the generators, $\{1\}\cup\{\xi_\b\st\b\le\a\}$, of $B_\a$ over $A$ so $H_\a=H|_{B_\a}$.\eop

\cor 5.8: Let $B$ be a standard normed extension of $A$. Then $K(A)=A\cap K(B)$.
\pf $$K(A)=\bigcap_{H\in\Omega(B)}\ker H|_A=\bigcap_{H\in\Omega(B)}A\cap\ker H=A\cap K(B).
\eqno\square$$

\vskip 10pt
\def\r{\rho}
\def\rs{\rho^*}
\noindent The second application concerns the similarity with Cole's construction in the case when $(A, X)$ is a
uniform algebra for we have (in the above notation and the notation of chapter 4) a map $\r\mapto{X^\GU}\OB{(\kappa,
\lambda)}{\amalg^{-1}(\eh \kappa, \lambda)}$. $\r$ is a homeomorphism onto its image in $\OB$ and induces a homomorphism
$\rs\map{C(\OB)}{C(X^\GU)}$ as in chapter 4. Moreover we have the corresponding

\lem 5.9: $\rs(\hat B)\subseteq A^\GU$.
\pf The proof is a straghtforward induction. Let $J=\{\b\in\GU\st\fa b {B_\b} \rs(\hat b)\in A^\GU\}$ and assume $[\ao, \a)\subseteq J$.

Let $b\in B_\a$. Let $n=n(\a)$, the degree of $\al$. Then there are $b_0,\ldots,b_{n-1}\in A_\a\st
b=b_0+\cdots+b_{n-1}\xi_\a^{n-1}$. So
$$\rs(\hat b)=\rs\left(\widehat{b_0}\right)+\cdots+\rs\left(\widehat{b_{n-1}}\right)
\rs\left( \widehat{\xi_\a} \right)^{n-1}$$
and by the inductive hypothesis it is enough to show that $\rs\left(\widehat{\xi_\a} \right)\in A^\GU$.
Now $\fa{(\kappa, \lambda)}{X^\GU}$
$$\rs\left( \widehat{\xi_\a} \right)(\kappa,\lambda)= \widehat{\xi_\a}\left( \amalg^{-1}(\eh\kappa,\lambda) \right)=
\amalg^{-1}(\eh\kappa,\lambda)(\xi_\a) =: H_\a(\xi_\a)$$
where (as in the proof of Theorem 5.6) $H_\a:=\amalg^{-1}(\eh\kappa,\lambda)|_{B_\a}=\omega\circ\phi_\a$ for $\omega\in\Omega\left(\Aa[x]/\prin\right)$ where in particular we have $\omega(\xb)=\lambda_\a$. Therefore $H_\a(\xi_\a)=\omega(\xb)=\lambda_\a=\pa(\kappa,\lambda)$. So (as expected) $\rs\left(\widehat{\xi_\a}\right)=\pa\in A^\GU$.

Thus $\a\in J$ so $J=\GU$.\eop

Moreover the following is commutative, because $\fa{f_0}A, (\kappa,\lambda)\in X^\GU$
$$\left(\rs\circ\widehat{f_0}\right)(\kappa, \lambda)=\amalg^{-1}(\eh\kappa,\lambda)(f_0)=\eh\kappa(f_0)=f_0(\kappa)=\pi^*(f_0)(\kappa,\lambda).$$

$$\matrix{B&\mapright\Gamma&\hat B &\mapright\rs&A^\GU\cr
	&\mapnw\iota&&\mapne{\pi^*}&\cr
	&&A&&\cr}$$

\cor 5.10: Let $(A, X)$ be a natural uniform algebra. Then $\rs$ is isometric and $A^\GU=\overline{B\char'136}$
\pf This follows from exactly the same reasoning as for simple extensions.\eop

\noindent As a last application of Theorem 5.6 we generalise in Theorem 5.12 Arens' and Hoffman's Theorem 2.10 from chapter 2.

\lem 5.11: Let $B$ be a standard extension of $A$ relative to $(B_\a)_{\a\in\GU}$ and $d\in A$ not be zero or a  zero divisor in $A$. Then $d$ is not a zero divisor in $B$.
\pf We use induction again but for some variation in the presentation suppose $d$ is a zero divisor in $B$ but not
in $A$. Recall that $B=\bigcup_{\a<\b}B_\a$. Let $\a$ be the least element of $\GU$ such that $dd\pri=0$ for some
$d\pri\in B_\a-\{0\}$.

As in the proof of Lemma 5.9 $\phi_\a\map{B_\a}{A_\a[x]/\prin}$ is an isomorphism over $\Aa$ so there are elements $ d_0,\ldots,d_{n(\a)-1} \in\Aa\st d\pri=d_0+\cdots+d_{n(\a)-1}\xi_\a^{n(\a)-1}$. Thus we have
$$0=\phi_\a(dd\pri)=dd_0+\cdots+dd_{n(\a)-1}\xb^{n(\a)-1}$$
By the uniqueness (discussed at the beginning of chapter 2) of this representation of the coset ($n(\a)$ is the degree of $\al$) we have $dd_j=0\;(j=0,\ldots,n(\a)-1)$. By assumption $d$ is not a zero divisor in $\Aa$ so $d_j=0\;(j=0,\ldots,n(\a)-1)$. So we have the contradiction $d\pri=0$.
\eop
 
\thm 5.12: Let $B$ be the standard normed extension of $A$ generated by $\GU$. Suppose that for every
$\a\in\GU\; d_\a:=\hbox{discr}(\al)$ is not zero or a zero divisor. Then $B$ is tractable if $A$ is.

\pf Let $A$ be tractable. Let $J=\{\b\in\GU\st B_\b\hbox{ is tractable}\}$ and suppose $[\ao,\b)\subseteq J$.

Let $a\in K(A_\b):=\cap M(A_\b)$. Now $\exists\,\a<\b\st a\in B_\a$ (or $a\in A$; it makes no difference in the following argument but where we write $B_\a$ we do mean `$B_\a$ or $A$').

Let $H_\a\in\Omega(B_\a)$.

Now $A_\b$ is a standard normed extension of $A$ relative to $(B_\gamma)_{\gamma<\b}$ so by Corollary 5.7, $$\te H{\Omega(A_\b)}\st H_\a=H|_{B_\a}.$$ \noindent
$\Rightarrow\;H_\a(a)=H(a)=0$. Thus $a\in K(B_\a)=(0)$ since $B_\a$ is tractable by hypothesis.

So $A_\b$ is tractable. Now $d_\b$ is not zero or a zero divisor in $A_\b$ by Lemma 5.11 so by Theorem 2.10, $A_\b[x]/(\beta(x))$, which is isometrically isomorphic to $B_\b$, is tractable.

By the principle of transfinite induction $J=\GU$ and it now follows from Corollary 5.7 again (with a proof as in the case of $A_\a$) that $B$ is tractable.\eop

\noindent In fact the converse is also true and the proof is another straightforward induction.
However it relies on the result being true in the simple case (Corollary 9.3 of [{\bf 7}]) so we omit this.

\vfill\eject

\def\GA{{\eulerfk U}}
\def\ao{\alpha_0}
\def\b{\beta}
\def\g{\gamma}
\def\u{\upsilon}
\def\k{\kappa}
\def\ps#1#2{\pi_{#1,#2}}
\def\pab{\pi_{\a,\b}}
\def\pss#1#2{\pi^*_{#1,#2}}
\def\psab{\pi^*_{\a,\b}}
\def\Hom{{\rm Hom}}
\def\Xfam{(X_\a)_{\a<\u}}
\def\Afam{(A_\a)_{\a<\u}}
\def\id{{\rm id}}

\def\fabu{\forall\,\a\le\b<\u\;}

\noindent {\bigtenbf Chapter 6}
\vskip 7pt
\noindent {\bigtenrm Integral Closure of Uniform Algebras}
\vskip 15mm

\noindent {\bf 1. Introduction}\vskip 10pt

\noindent In Cole's thesis the main example features a method for constructing uniform algebra extensions closed under
taking square roots and Lindberg remarks in [{\bf 11}] that the same techniques can be used to construct an integrally
closed uniform algebra which isometrically extends any given uniform algebra, $(A, X)$. We give the details explicitly here.
\vskip 15mm

\noindent {\bf 2. Inverse and Direct Systems}\vskip 10pt

\noindent The basic result is that there is a canonical way of forming transfinite extensions of uniform algebras
and it is this pre-existing construction which Cole exploits in [{\bf 10}]. In anticipation of the definition to be given below, the result is that `direct limits' of uniform algebras always exist (we shall prove a weaker result here). For the general categorical definitions of direct and inverse systems and their limits we refer the reader to [{\bf 15}], but it is easy to see what they are from the definitions below.

\vskip 15pt

\hrule\vskip5pt\noindent
Throughout the rest of this chapter, $\u$ will denote a fixed ordinal number greater than zero.
\vskip5pt\hrule\vskip 20pt

\dfn 6.1: Let $\Xfam$ be a family of topological spaces and $(\pab)_{\a\le\b<\u}$ a family of continuous maps where $\forall\,\a\le\b<\u\;\pab\in C(X_\b,X_\a)$. We say that $\left(\Xfam, (\pab)_{\a\le\b<\u}\right)$ is an {\it inverse system} of topological spaces if
\item{(i)} $\forall\,\a<\u\ \ps\a\a=\id_{X_\a}$
\item{(ii)} $\forall\,\a\le\b\le\g<\u\ \ps\a\g=\pab\circ\ps\b\g$.

An {\it inverse limit} for the system above is a pair $\left(X_\u, (\ps\a\u)_{\a<\u}\right)$ where $X_\u$ is a topological space and $\forall\,\a<\u\;\ps\a\u\in C(X_\u,X_\a)$ are such that
$$\forall\,\a\le\b<\u\qquad \ps\a\u=\pab\circ\ps\b\u$$ and such that the following universal property is satisfied:
Whenever $Y_\u$ is a topological space and $\exists_{\a<\u}\,\rho_{\a,\u}\in C(Y_\u, X_\a)\st \forall\,\a\le\b<\u\ \rho_
{\a,\u}=\pab\circ\rho_{\b,\u}$ then there is a unique $\rho_\u\in C(Y_\u, X_\u)\st$
$$\forall\,\a<\u\qquad\ps\a\u\circ\rho_\u=\rho_{\a,\u}$$
This is written $$X_\u=\il\a\, X_\a.$$

\noindent In general there is no guarantee that an inverse limit will exist or be non-trivial. However we have the important result below:

\thm 6.2: Let $\left(\Xfam, (\pab)_{\a\le\b<\u}\right)$ be an inverse system of non-empty compact Hausdorff spaces. Then an inverse limit exists and is a non-empty compact Hausdorff space.
\pf (From [{\bf 16}].): Define $K=\times_{\a<\u}X_\a$ and let $\pi_\a\map K {X_\a}$ be the canonical projection ($\a<\u$). Set
$$\eqalign{X_\u&:=\left\{\kappa\in K\st\forall\,\a\le\b<\u\;(\pab\circ\pi_\b)(\kappa)=\pi_\a(\kappa)\right\}\cr
\ps\a\u&:=\pi_\a|_{X_\u}\qquad(\a<\u).\cr}$$
We show that $\left(X_\u, (\ps\a\u)_{\a<\u}\right)$ is an inverse limit for the system.

A product of Hausdorff spaces is Hausdorff so $K$ is Hausdorff and $K$ is compact by Tychonoff's theorem. We have that $X_\u=\bigcap_{\b<\u}F_\b$ where $\forall\,\b<\u\;F_\b=\left\{\kappa\in K\st\forall\,\a\le\b\;(\pab\circ\pi_\b)(\kappa)=\pi_\a(\kappa)\right\}$. By exercise 3C of [{\bf 17}] the set on which two continuous maps into a Hausdorff space agree is closed so the $F_\b$ are closed. Hence $X_\u$ is a closed and therefore compact subset of $K$.

To check that $X_\u$ is non-empty it is enough, since it is now known to be compact and therefore has the finite intersection property, to show that $\forall\,\b_1,\ldots,\b_n<\u\;\bigcap_{k=1}^n F_{\b_k}\ne\emptyset$.

We first need to know that $\forall\,\b<\u\;F_\b\ne\emptyset$.

Let $\b<\u$ and choose $\k_\b\in X_\b$ and for $\a\le\b$ let $\k_\a=\pab(\k_\b)$. Let $\k_\a$ be arbitrary in $X_\a$ for $\b<\a<\u$. Then $\k\in K$ and $\forall\,\a\le\b\;(\pab\circ\pi_\b)(\k)=\pab(\k_\b)=\k_\a=\pi_\a(\k)$ so $\k\in F_\b$.

Now let $\b=\max_{k=1,\ldots,n}(\b_k)$. Let $\k\in F_\b$. We have
$$\eqalign{&\qquad\forall\,\a\le\b\qquad (\pab\circ\pi_\b)(\k)=\pi_\a(\k)\cr
&\hbox{and so}\qquad\fa j{\{1,\ldots,n\}},\;\a\le\b_j\qquad\left(\ps\a{\b_j}\circ\pi_{\b_j}\right)(\k)=\left(\ps\a{\b_j}\circ\ps{\b_j}\b\circ\pi_\b\right)(\k)\cr
&\hbox{whence}\qquad\fa j{\{1,\ldots,n\}},\;\a\le\b_j\qquad\left(\ps\a{\b_j}\circ\pi_{\b_j}\right)(\k)=(\pab\circ\pi_\b)(\k)=\pi_\a(\k)\cr
&\hbox{and so}\qquad \k\in\bigcap_{j=1}^n F_{\b_j}.\cr}$$
Thus $\emptyset\ne\bigcap_{j=1}^n F_{\b_j}$.

It remains to check that $\left(X_\u, (\ps\a\u)_{\a<\u}\right)$ satisfies the rest of the conditions in 6.1 above. The maps are clearly continuous by the definition of the product topology. Now suppose that $\left(Y_\u, (\rho_{\a,\u})_{\a<\u}\right)$ is a topological space and family of continuous maps with $\forall\,\a\le\b<\u\ \rho_{\a,\u}=\pab\circ\rho_{\b,\u}$. Define
$$\rho_\u\mapto{Y_\u}{X_\u}y{\left(\rho_{\a,\u}(y)\right)_{\a<\u}}$$
Then $\fabu (\pab\circ\pi_\b)(\rho_\u(y))=\pab(\rho_{\b,\u}(y))=\rho_{\a,\u}(y)=\pi_\a(\rho_\u(y))$ so $\rho_\u$ is well-defined and has $\forall\,\a<\u\quad\ps\a\u\circ\rho_\u=\rho_{\a,\u}$. $\rho_\u$ is continuous since $\forall\a<\u\ \ps\a\u\circ\rho_\u$ is continuous and $K$ has the initial topology induced by $(\pi_\a)_{\a<\u}$ so $X_\u$ has initial topology induced by $\left(\pi_\a|_{X_\a}\right)_{\a<\u}=(\ps\a\u)_{\a<\u}$. Uniqueness of $\rho_\u$ is clear.\eop

\def\thab{\theta_{\a,\b}}
\def\ths#1#2{\theta_{#1,#2}}

\def\phs#1#2{\phi_{#1,#2}}
\def\thu{\theta_\u}
\def\thau{\theta_{\a,\u}}
\dfn 6.3: Let $\left(\Afam, (\thab)_{\a\le\b<\u}\right)$ be a family of algebras and homomorphisms where $\fabu\thab\map{A_\a}{A_\b}$. This is a {\it direct system} if
\item{(i)} $\forall\,\a<\u\qquad \ths\a\a=\id_{A_\a}$
\item{(ii)} $\forall\,\a\le\b\le\g<\u\qquad \ths\a\g=\ths\b\g\circ\thab$.

A {\it direct limit} of this system is $\left(A_\u, (\ths\a\u)_{\a<\u}\right)$ for an algebra $A_\u$ and homomorphisms $\ths\a\u\map{A_\a}{A_\u}$ such that
$$\fabu\qquad \ths\a\u=\ths\b\u\circ\thab$$ and whenever $B_\u$ is an algebra with homomorphisms $\phs\a\u\map{A_\a}{B_\u}\ (\a<\u)$ satisfying
$$\forall\,\a\le\b<\g\qquad \phs\a\u=\phs\b\u\circ\ths\a\b$$
then
$$\exists\vert\,\thu\in\Hom(A_\u,B_\u)\st\qquad\forall\,\a<\u\quad \phs\a\u=\thu\circ\thau.$$
We write $$A_\u=\dl\a\,A_\a.$$

\noindent It is clear what these definitions in other specific categories of algebras should be. In fact for our purposes it will be
sufficient to restrict attention to the category, {\cal U}, of uniform algebras and continuous unital homomorphisms between them of norms at most 1 (i.e. multiplicative contractions).

\vskip 10pt\noindent
Now let $\left((A_\a, X_\a)_{\a<\u}, (\pab)_{\a\le\b<\u}\right)$ be a family of uniform algebras such that $\left(\Xfam, (\pab)_{\a\le\b<\u}\right)$ is an inverse system. Then
$$\eqalign{\fabu\qquad\pab\hbox{ induces }&\psab\in\Hom_\GA(A_\a,A_\b)\quad\hbox{ defined by}\cr
&\psab\mapto{A_\a}{A_\b}f{f\circ\pab}.\cr}$$
The success of Cole's method rests (among other things) on:

\thm 6.4: Let $\left((A_\a, X_\a)_{\a<\u}, (\pab)_{\a\le\b<\u}\right)$ be a family of uniform algebras and suppose $$\left(X_\u, (\ps\a\u)_{\a<\u}\right)$$ is an inverse limit
of the inverse system $\left(\Xfam, (\pab)_{\a\le\b<\u}\right)$ where all the $\pab$ are surjective $(\a\le\b<\u)$. Let $A_\u$ be the closed subalgebra of $C(X_u)$ generated by $I_0:=\bigcup_{\a<\u}\pss\a\u(A_\a)$. Then
$\left(A_\u, (\pss\a\u)_{\a<\u}\right)$ is a direct limit of the dual system $\left(\Afam, (\psab)_{\a\le\b<\u}\right)$ (in {\cal U}) and $(A_\u,X_\u)$ is a uniform algebra.

\pf (Adapted from the sketch in [{\bf 10}].) It is immediate that $(A_\u,X_\u)$ is a uniform algebra; $X_\u$ is a
non-empty compact Hausdorff space by Theorem 6.2 and $A_\u$ is already a uniformly closed algebra.
It contains the constants because $\pss0\u(1)=1\in A_\u$. $A_\u$ separates the points of $X_\u$: let
$\k,\k^\prime\in X_\u$ with $\k\ne\k^\prime$. So there exists $\a<\u\st\k_\a\ne\k_\a\pri$. But $A_\a$ is a
uniform algebra so for some $f\in A_\a$, $f(\k_\a)\ne f(\k\pri_\a)$. We have $\pss\a\u(f)\in A_\u$ and
$\pss\a\u(f)(\k)=(f\circ\ps\a\u)(\k)=f(\k_\a)\ne\pss\a\u(f)(\k\pri)$.

We now verify that $\left(A_\u, (\pss\a\u)_{\a<\u}\right)$ is a direct limit of the dual system in {\cal U}. By Lemma 9 of
[{\bf 16}] (p. 210), whose proof we do not reproduce here, the surjectivity of the maps $(\pab)_{\a\le\b<\u}$ implies that
$\ps\a\u$ is surjective $(\a<\u)$. Hence $\pss\a\u$ is isometric $(\a<\u)$. So for $\a<\u\;\norm{\pss\a\u}=1$. Also
$$\fabu\qquad \pss\b\u\circ\psab=(\pab\circ\ps\b\u)^*=\pss\a\u.$$
It therefore remains to check that $\left(A_\u, (\pss\a\u)_{\a<\u}\right)$ solves the universal problem in definition 6.3.

Let $\left( (B_\u,Y_\u), (\phs\a\u)_{\a<\u}\right)$ be another pair of uniform algebras and family of homomorphic contractions, $\phs\a\u\map{A_\a}{B_\u}\;(\a<\u)$, with
$$\fabu\qquad \phs\a\u=\phs\b\u\circ\psab.$$	\def\pus{\pi_\u^*}
We must show $$\exists\vert\,\pus\in\Hom(A_\u, B_\u)\st\forall\,\a<\u\quad\phs\a\u=\pus\circ\pss\a\u.\eqno(\dag)$$

The uniqueness part is clear because $I_0$ generates $A_\u$. Also notice that $I_0$ is a subspace of $A_\u$: let $g_1,g_2\in I_0, \lambda\in\CN$.

There exist $\a,\b<\u\st g_1\in\pss\a\u(A_\a), g_2\in\pss\b\u(A_\b)$.

So there are $f_1\in A_\a,f_2\in A_\b\st g_1=\pss\a\u(f_1),g_2=\pss\b\u(f_2)$.

We have $\lambda g_1=\pss\a\u(\lambda f_1)\in I_0$. Now let $\a\le\b$ without loss of generality. So $\pss\a\u=\pss\b\u\circ\psab$ whence $g_1+g_2=\pss\b\u\left(\psab(f_1)+f_2\right)\in I_0$.

So if we can define $\pus$ on $I_0$ to be a homomorphism with norm at most 1 satisfying (\dag) then (since $B_\u$ is a
Banach space) there is a unique continuous linear extension to a map $A_\u\to B_\u$ with norm at most 1. By continuity the extension will be a homomorphism.

Set $T_0\mapto{I_0}{B_\u}g{\phs\a\u(f)}$ where $g=\pss\a\u(f), f\in A_\a, \a<\u$. $T_0$ is well-defined: suppose $\a\le\b<\u$ and $g=\pss\a\u(f)=\pss\b\u(f\pri)$ for $f\in A_\a,f\pri\in A_\b$. Now $\pss\a\u=\pss\b\u\circ\psab$ so $\pss\a\u(f)=\pss\b\u(\psab(f))=\pss\b\u(f\pri)$. But $\ps\b\u$ is surjective so $\pss\b\u$ is a monomorphism (in fact isometric, as noted above). So $f\pri=\psab(f)$\par\noindent
$\Rightarrow\qquad\phs\a\u(f)=(\phs\b\u\circ\psab)(f)=\phs\b\u(f\pri)$.

Similar calculations show that $T_0$ is a homomorphism. We also have that $T_0$ is a contraction for if $g=\pss\a\u(f), f\in A_\a$ then $$\norm{T_0(g)}=\norm{\phs\a\u(f)}\le\norm f=\norm{\pss\a\u(f)}=\norm g.$$
So the required map $\pus\map{A_\u}{B_\u}$ exists.\eop
\vskip 15mm

\noindent {\bf 3. Universal Root Algebras}\vskip 10pt

\noindent We can now apply the machinery above to obtain a `universal root algebra' (an integrally closed extension algebra)
by mimicking Cole's work.
 This is almost certainly repeating the work of Lindberg.

\def\w{\omega_1}
\thm 6.5 Let $(A,X)$ be a uniform algebra. Then there exists an integrally closed uniform algebra which is an isometric extension of $A$.

\pf Let $\w$ denote the first uncountable ordinal. Now let $\g<\w$ and suppose we have chosen uniform algebras, $(A_\a, X_\a)_{\a<\g}$, and continuous surjections, $(\pab)_{\a\le\b<\g}$, such that
\item{(i)} $\left((X_\a)_{\a<\g}, (\pab)_{\a\le\b<\g}\right)$ is an inverse system
\item{(ii)} $(A_0, X_0)=(A, X)$
\item{(iii)} $\forall\,\b<\g$, if $\exists\,\a<\b\st\b=\a+1$ \footnote{*}{See appendix 3} then $(A_\b, X_\b)$ is the
Cole-extension as in Theorem 4.1 taken with respect to $\GU_\a$, the set of all monic polynomials over
$A_\a$ and $\forall\,\a\le\b<\g\;\pab$ is the associated projection $(X_\a)^{\GU_\a}=X_\a\times\CN^{\GU_\a}\to X_\a$.
\item{(iv)} $\forall\,\b<\g$, if $\b$ is a limit ordinal then $A_\b$ is the direct limit of $(A_\a, X_\a)_{\a<\b}$
and $(\ps\a\u)_{\a<\u}$ are as described in Theorem 6.4.

There are two cases. If $\g=\b+1$ for some $\b<\g$ then we (are forced to) define $(A_\g, X_\g)$ to be the uniform
algebra generated according to Theorem 4.1 by the set of all monic polynomials over $A_\b$. For $\a\le\g$ set
$$\ps\a\g=\cases{\id_{X_\g} & $\a=\g$\cr
	\hbox{the coordinate projection }\ps\b\g\map{X_\g}{X_\b} & $\a=\b$\cr
	\pab\circ\ps\b\g&$\a<\b$.\cr}$$
All maps on the right hand side are surjective and continuous so the same is true of $\ps\a\g\;(\a\le\g)$.

The other case is when $\g$ is a non-zero limit ordinal. Define $\left((A_\g, X_\g), (\ps\a\g)_{\a<\g}\right)$ to be
the uniform algebra constructed out of the existing system according to Theorem 6.4 (and let $\ps\g\g=\id_{X_\g}$, the
identity map on $X_\g$).

In each case $\left((X_\a)_{\a\le\g}, (\pab)_{\a\le\b\le\g}\right)$ is an inverse system of compact Hausdorff spaces and continuous surjections and (i)-(iv) are plainly satisfied by the extended system of uniform algebras. By the transfinite recursion theorem we obtain such a family of uniform algebras with $\g=\w$.

It is clear that $\w$ is a limit ordinal. Therefore by (iv),
$$A_{\omega_1}=\overline{\bigcup_{\a<\w}\pss\a{\omega_1}(A_\a)}$$
the closure being taken in $C(X_{\omega_1})$. This follows from the facts that $\bigcup_{\a<\w}\pss\a{\omega_1}(A_\a)$ is a
subalgebra (similar to the result seen in the proof of Theorem 6.4) and that the closure of a subalgebra is a subalgebra.
We also saw in the proof of Theorem 6.4 that $\pss\a{\omega_1}$ is an isometry ($\a<\w$). Hence $\pss\a{\omega_1}(A_\a)$ is closed $\forall\,\a<\w$ (since $A_\a, C(X_{\omega_1})$ are complete).

Also the relation $\forall\,\a\le\b<\w\ \pss\a{\omega_1}=\pss\b{\omega_1}\circ\psab$ shows that
$\forall\,\a\le\b<\w\ \pss\a{\omega_1}(A_\a)\subseteq\pss\b{\omega_1}(A_\b)$. Hence by Theorem A3.20, $\bigcup_{\a<\w}\pss\a{\omega_1}(A_\a)$ is closed. Therefore
$$A_{\omega_1}=\bigcup_{\a<\w}\pss\a{\omega_1}(A_\a).$$

Now if $\mu(x)$ is a monic polynomial over $A_{\omega_1}$ then it is monic over $\pss\a{\omega_1}(A_\a)$ for some $\a<\w$ and therefore has a root in $\pss\b{\omega_1}(A_\b)$ where $\b=\a+1$. So $A_{\omega_1}$ is integrally closed. Finally note that $\pss0{\omega_1}$ gives an isometric embedding $A\hookrightarrow A_{\omega_1}$. \eop
\vskip 10pt
We therefore have refined Theorem 5.5 to the category of uniform algebras. One can also investigate the inheritance
properties of these `universal root extensions'. For example Theorem 4.2 generalises: if $A$ is natural on $X$ then
$A_{\omega_1}$ is natural on $X_{\omega_1}$. This can be proved in exactly the same way as the corresponding
result for the square-root closed extensions mentioned before; see for example [{\bf 18}] p. 103.

\vfill\eject

\noindent {\bigtenbf Chapter 7}
\vskip 7pt \noindent {\bigtenrm Conclusion}
\vskip 10mm

\noindent {\bf 1. Conclusions}\vskip 10pt

\noindent We have explored two methods for extending uniform algebras by roots of monic polynomials, one of which applies to general commutative normed unital algebras. We have seen that integrally closed extensions always exist in both categories of algebras.

The two constructions have well-understood properties, and, as we have seen lead to essentially different solutions.

We have taken a simplified view of the literature. There have been many contributions (mainly by Lindberg) on properties of Banach algebras with respect to Arens-Hoffman extensions which we have not considered. For example (for appropriate initial algebras) naturality, normality, regularity, the \v Silov boundary, indecomposability and involutions. (See [{\bf 7}].)

Another topic not covered here is on `integral extensions'. Standard extensions are examples of integral extensions and they can be used to investigate the purely algebraic questions on them; see [{\bf 11}].

\vskip 15mm

\noindent {\bf 2. Context}\vskip 10pt

\noindent Apart from the minor application in Corollary 3.16, Arens-Hoffman extensions' mainly feature in the literature in connection with Galois theory for Banach algebras where they play a very important role. However in some versions of this they are not sufficient (there are `Galois extensions which are not Arens-Hoffman extensions' mentioned in [{\bf 19}]).

Arens-Hoffman extensions also feature in some examples in the literature on function algebras. For example, Karahanjan uses them
in [{\bf 14}] to construct an example of a uniform algebra with certain properties.
This was the area of Cole's work in [{\bf 10}] where he used a universal root algebra to disprove a certain conjecture.

\vskip 15mm

\noindent {\bf 3. Further Questions}\vskip 10pt

\noindent In addition to the properties mentioned in the first section there are other inheritance properties of Arens-Hoffman extensions it might be interesting to investigate; for example the popular $C$* property. Also the work in section 5.4 suggests that many of these results will lift to standard extensions.

The corresponding results for extension of uniform algebras by Cole's methods have not been studied as much in the literature so
it is important to know precisely when such extensions can be realised via Arens-Hoffman extensions.

We would also like to be able (if possible) to generalise the result at the end of chapter two so that every tractable normed algebra can be embedded in a {\it tractable} normed algebra extension which contains roots of an arbitrary large set of monic polynomials (which may have to be restricted in some other way).

\vfill\eject

\noindent {\bigtenbf Appendix 1}
\vskip 7pt \noindent {\bigtenrm Assumed Gelfand Theory}
\vskip 15mm

\noindent Gelfand theory is a tool used to investigate normed algebras.
It identifies arbitrary normed algebras with complex function algebras via
a map, $\Gamma$, to be defined below. We list the definitions and results which
are used in the dissertation. This material is standard and proofs can be found
in many introductory books on functional analysis, for example [{\bf 5}] chapters 12 and 13 or [{\bf 20}] chapters 1 and 2.
\par
Let $A$ be a commutative unital Banach algebra.

\dfn A1.1: The {\it maximal ideal space} and {\it character (carrier) space}
of $A$ are respectively
$$\eqalign{M(A)&:=\{I\subseteq A\colon I \hbox{ is a maximal ideal of }A\},\qquad\hbox{ and}\cr
\Omega(A)&:=\hbox{Hom}_\CN(A,\CN)-\{0\}.\cr}$$
The topologies on these two sets are defined below. Recall that
$A^*={\cal B}(A, \CN)$ has a topology called the weak *-topology, the initial
topology induced by the evaluation-at-$a$ linear functionals ($a\in A$).

\thm A1.2: $\Omega(A)\subseteq A^*$.\eop

\thm A1.3: $M(A)\ne\emptyset$ and $\mu\colon\Omega(A)\to M(A);\,\omega
\mapsto\ker\omega$ is a bijection.\eop

\noindent The topology on $\Omega(A)$ is the subspace topology relative to the weak
*-topology. The topology on $M(A)$ is the one obtained by using $\mu$ to
identify the two sets.

\thm A1.4: $\OA$ is a compact Hausdorff subspace of $A^*$ in the weak
*-topology.\eop

\dfn A1.5: The {\it Gelfand transform} of $A$ is the map $\Gamma\colon
A\to C(\OA);\,a\mapsto \ha$ where
$$\forall\,a\in A,\omega\in\OA\qquad \ha(\omega):=\omega(a).$$

\thm A1.6: \item{(i)} $\Gamma$ is a contractive (norm-decreasing) homomorphism
\item{(ii)} $\forall\,a\in A\ \norm \ha = r(a) = \lim_{n\to+\infty}
		\norm{a^n}^{1\over n}$
\item{(iii)} $\forall\,a\in A\ \im\ha=\sigma(a)$
\item{(iv)} $\ker\Gamma=\cap M(A).$\eop

\dfn A1.7: $A$ is {\it semisimple} if $\cap M(A)=\{0\}$.

\cor A1.8: \item{(i)} $\Gamma$ is injective $\iff$ $A$ is semisimple
\item{(ii)} $\Gamma$ is isometric $\iff\ \forall\,a\in A\ \norm
		{a^2}=\norm a^2$.\eop

\vfill\eject

\def\sb{{\bf s}}
\def\sbp{\sb^\prime}
\def\ab{{\bf a}}
\def\bb{{\bf b}}

\noindent {\bigtenbf Appendix 2}
\vskip 7pt \noindent {\bigtenrm Gelfand Theory for Normed Algebras}
\vskip 15mm

\noindent While Gelfand theory is usually expounded for Banach algebras it is
interesting to study its generalisation to normed algebras in its own right
as well for its use as a tool to investigate incomplete normed algebras. Most of
the results presented below seem to be well-known but not to be found explicitly
stated in the literature. They are obtained by modifying the standard proofs for the
complete case or by considering the algebra's completion. Accordingly we begin
by recalling the definition of this.
\dfn A2.1: Let $A$ be a normed algebra. Then a  {\it completion} of $A$ is a pair
$(\tilde A, i)$ where $\Ac$ is a Banach algebra and $i$ is an isometric
monomorphism $A\to\Ac$ with dense image.\par
There is a well known construction for normed spaces. One defines $\Ac:=S/\sim$ where $S=\{\sb\in A^\NN\,\colon\,
\sb\hbox{ is a Cauchy sequence}\,\}$ and $\sim$ is the equivalence relation on $S$ given
by $$\sb\sim \sb^\prime\qquad :\Leftrightarrow\qquad \lim_{n\to+\infty}
\Vert s_n-s_n^\prime\Vert=0\qquad\qquad{\rm where}\ \sb=(s_n), \sbp=(s_n^\prime)$$
It is standard that completions are unique up to isometric isomorphism. By this we mean that if $(B, j)$ is another completion then there is an isometric isomorphism $\phi\map {\tilde A} B$ such that the following is commutative:
$$\matrix{ \tilde A&\mapright\phi&B\cr
	\mapup i&\mapne j&\cr
	A&&\cr}$$

\noindent 
For $\ab,\bb\in S\ (\ab=(a_n), \bb=(b_n)),\lambda\in \CN$ the normed space 
structure is given by
$$\eqalign{[(a_n)]+[(b_n)]&:=[(a_n+b_n)]\cr
\lambda[(a_n)]&:=[(\lambda a_n)]\cr
\Vert [(a_n)]\Vert&:={\rm lim}_{n\to+\infty}\Vert a_n\Vert\cr}$$
It is a standard exercise to check that

$$\eqalign{[(a_n)][(b_n)]&:=[(a_nb_n)]\cr}$$
properly defines a multiplication on $\Ac$ in such a way that $i$ is an
algebra homomorphism. It is also clear that normed algebra completions are
unique to isometric isomorphism. They are clearly commutative if $A$ is
commutative. If $A$ is unital then so is $\Ac$ and $i$ is too
(i.e. $1_{\Ac}=i(1_A)$).

\vskip 15pt
\hrule\vskip5pt\noindent
In the rest of this appendix $A$ will denote a complex normed algebra, not necessarily commutative or unital. \vskip5pt\hrule\vskip 20pt

\noindent Although the applications in the dissertation are only to unital algebras it is not much more difficult to develop the Gelfand theory for non-unital normed algebras. There is an existing theory for non-unital Banach algebras (see for example [{\bf 20}] chapter 1) in which the maximal ideals are replaced by the modular maximal ideals.

\dfn A2.2: Let $R$ be a commutative algebra and $J$ be a maximal ideal of $R$. $J$ is called {\it modular} if for some
$u\in J$ we have $\fa r R\ r-ru\in J$.\vskip 10pt

\noindent Clearly all ideals are modular if $R$ has a unity.

\noindent
When $A$ is not complete it is no longer automatic that a complex homomorphism
is continuous. We therefore only consider the continuous characters on $A$:
$$\Omega(A):=\{\omega\in{\rm Hom}_\CN(A, \CN)\colon \omega
\ \hbox{is continuous,}\ \omega\ne 0\}$$
In this way $\OA\subseteq A^*$ and we can again exploit the topology on $A^*$.
\lem A2.3 (Gelfand-Mazur): Let $A$ be a normed, unital division algebra. Then
there exists an isometric isomorphism $A\to\CN$.
\pf This will be as for Banach algebras if we can show that again in normed algebras
every $a\in A$ has non-empty spectrum, $\sigma(a)$. This is the content of the Gelfand-Mazur
theorem for Banach algebras. It can be shown directly for normed algebras (for example
in [{\bf 2}] p. 22 Theorem 5.6) but the following argument, pointed out to me by Dr. Feinstein,
avoids this. By the Banach algebra result $\sigma_{\Ac}(a)\ne\emptyset$ where $\Ac$ is
as above. But since there is a unital monomorphism $A\hookrightarrow\Ac$ we have
$\sigma_{\Ac}(a)\subseteq\sigma_A(a)$ so $\sigma(a)\ne\emptyset$ as required.\eop

\par\noindent One comment is in order: if the normed algebra has a unity but is not unital (i.e.
$1<\Vert 1_A\Vert$) then the topological isomorphism of the theorem still exists but 
may not be isometric. (The map is $A\to\CN;a\mapsto\lambda_a$ where $\lambda_a$ 
is the unique $\lambda\in\CN\colon a=\lambda 1_A$.)

\lem A2.4: Let $A$ be a commutative normed algebra and let 
$$M(A):=\{I\subseteq A\colon I
\hbox{ is a closed, modular maximal ideal of }A\}.$$ Then $M(A)\ne\emptyset\Leftrightarrow\Omega(A)\ne\emptyset$
and in this case $\mu\colon\Omega(A)\to M(A);\,\omega\mapsto{\rm ker}\,\omega$
is a bijection.
\pf First suppose $\Omega(A)\ne\emptyset$. Let $\omega\in\Omega(A)$. Since
$\omega\ne 0$ its kernel is a proper ideal of $A$. Suppose $J$ is an ideal of $A$
satisfying $\ker\omega\subset J\subset A$. Now $\omega$ is a ring epimorphism
$A\to\CN$ so $\omega(J)$ is an ideal of $\CN$. But $\CN$ is simple so
$\omega(J)=\left\{0\right\}$ or $\CN$.

If $\omega(J)=(0)$ then $J\subseteq\omega^{-1}(\omega(J))=\omega^{-1}((0))=\ker\omega$ which implies that
$\ker\omega=J$, a contradiction. The other case can't occur because if $\omega(J)=\CN$ then $0=\CN/\omega(J)\cong
A/J$ (since $\omega$ is an epimorphism). But $J\ne A$ so $A/J\ne 0$, a contradiction.

Thus $\ker\omega$ is maximal.

\par $\ker\omega$ is closed by the continuity of $\omega$.\par
Let $a\in A-\ker\omega$ and put $u=a/\omega(a)$. Then we have that
$$\forall\,b\in A\qquad b-ub\in\ker\omega$$ and so $\ker\omega$ is modular.
Thus $\ker\omega\in M(A)$. Note that this also shows that $\mu$ is well-defined and exists
whenever $\Omega(A)\ne\emptyset$. Also, the result that $\mu$ is always injective when
it exists does not require the completeness of $A$ so we omit the detail
(see for example [{\bf 5}] p. 321).\par
Suppose now that $M\in M(A)$. Since $A$ is commutative, $A/M$ is a field; the unity is
$M+u$ where $u\in A\colon\forall\,a\in A\ a-ua\in M$. By Lemma A2.3 there is a
topological isomorphism $\phi\colon A/M\to\CN$. The map $\omega:=\phi\circ\nu$ where
$\nu$ is the natural epimorphism is therefore a continuous character of $A$. (It is
standard that the natural epimorphism of a normed space onto its quotient by a closed
subspace is continuous.)\par
We have now shown $M(A)\ne\emptyset\Leftrightarrow\Omega(A)\ne\emptyset$.
It remains to check that when either set is non-empty, $\mu$ is surjective. But for $M\in M(A)$, setting $\omega=\phi\circ\nu$ as described above gives a continuous character on $A$ with $M\subseteq\ker\omega$. Since $\ker\omega\subset A$ and $M$ is maximal we have $M=\ker\omega=\mu(\omega)$ as required.\eop

\def\oc{\tilde \omega}
\lem A2.5: Let $A$ be a normed algebra with $\OA$ non-empty.
Every continuous character on $A$ extends to a character on $\Ac$.
\pf Let $\omega\in \OA$. We must show that $\exists\,\oc\in\Ac\,\colon\,\oc\circ i=\omega$
where $(\tilde A, i)$ is as above. But $\im i$ is dense in $\Ac$, $\omega\circ i^{-1}$ is a
bounded linear map $\im i\to\CN$ with norm $\Vert \omega\Vert$ (where $i^{-1}$ is the (induced)
isometric isomorphism $\im i\to A$). $\CN$ is complete so by the standard continuous
linear extension theorem, $\exists|\,\oc\in{\cal B}(\Ac, \CN)\colon\ \oc|_{\im i}=\omega$.
This continuous extension is clearly also a unital homomorphism satisfying 
$\oc\circ i=\omega$.\eop

As for Banach algebras we define, when $\OA\ne\emptyset$, for $a\in A$ the map
$$\ha :\OA\to\CN;\,\omega\longmapsto \omega(a).$$ Let the set of these maps be $\hA$.

\def\ti{\tau_i}
\lem A2.6: Let $A$ be a normed algebra with $\OA$ non-empty.
Let $\ti$ be the initial topology on $\OA$ induced by $\hA$. Then
$\OA$ is a locally compact Hausdorff space and is compact if $A$ is unital.
\pf First we observe that $\ti=\tau^*_{\OA}$ where $\tau^*$ is the weak *-topology
on $A^*$. This is because $\ti$ is the initial topology induced by the restrictions to
$\OA$ of the set of maps inducing the topology $\tau^*$ (see [{\bf 5}] exercise 19-C).
\par So it is immediate that $\OA$ is Hausdorff, since
$\tau^*$ is known to be. Now let $\omega\in\OA$.  \par By Lemma A2.5, $\exists\,\oc\in\Omega(\Ac)\colon \oc\circ i$.
By the result for Banach algebras, $\oc\in\Omega(\Ac)\;\Rightarrow\;\Vert\oc\Vert\le 1$. So
$\Vert \omega\Vert=\Vert\oc\circ i\Vert\le\Vert\oc\Vert\Vert i\Vert\le 1$. Therefore
$\OA\subseteq S^*:={\rm B}_{A^*}[0, 1]$. By the Banach-Alaoglu theorem, $S^*$ is compact
with respect to $\tau^*_{\OA}=\ti$. The rest of the proof (that $\OA\cup\{0\}$ is closed
in $S^*$ so that $\OA$ is locally compact) is exactly the same as for Banach algebras 
(see for example [{\bf 20}] \res1.3.5) so we omit this.\eop

\cor A2.7: $\fa\omega\OA\ \norm\omega=1$.\eop

\lem A2.8: Let $A$ be a normed algebra with $\OA$ non-empty. Then $\hA$ is a subalgebra
of $C_0(\OA)$ which separates the points of $\OA$ strongly.
\pf Let $\Gamma\colon A\to C_0(\OA);\,a\mapsto\ha$. It is clear that $\hA=\im\Gamma$ separates the points of
$\OA$ strongly: if $\omega_1,\omega_2\in\OA,\ \omega_1\ne\omega_2\Rightarrow\ \exists
a\in A\,:\,\omega_1(a)=\ha(\omega_1)\ne\ha(\omega_2)=\omega_2(a)$. That $\hA$ vanishes at
no point of $\OA$ is built into the definition of characters.\par
The verifications that $\Gamma$ is a homomorphism and that $\forall\,a\in A\ \ha$ is 
continuous and vanishes at infinity are the same as for
Banach algebras so we omit the detail. Since $\hA=\im\Gamma$, it is automatically
an algebra.\eop

\noindent To summarise our progress: we have established a norm decreasing homomorphism, $\Gamma$, from the normed algebra $A$ into a $C_0(X)$ for some locally compact Hausdorff space,
and that in fact it is given in the same way as for Banach algebras but restricting
attention to the continuous characters of $A$. We shall still call $\Gamma$ the Gelfand
transform. However it is now longer true that $\Gamma$ is injective iff $A$ is semisimple
and the relevant condition on $A$ is given by the following definition, given in 
[{\bf 1}]:
\dfn A2.9 ([{\bf 1}]): Let $A$ be a normed unital algebra. $A$ is called {\it tractable} if the intersection of its closed maximal ideals is $\{0\}$.\par
(Recall that the existence of the unit implies maximal ideals exist.) In fact by lemma
A2.4 we can generalise (for commutative $A$ at least) this to
\dfn A2.10 Let $A$ be a normed algebra. A is {\it tractable} if the intersection of its
closed modular maximal ideals, $M(A)$, is $\{0\}$.
\par There is a definition of modular ideals for non-commutative normed algebras but we shall not refer to it.
(The definition implies that $M(A)$ is non-empty (unless $A=0$) since the empty intersection of subsets of $A$
is $A$ itself.) It is consistent with the unital case because then every maximal ideal is modular with `modular unit' $1_A$. We have

\thm A2.11: Let $A$ be a normed algebra for which $\OA\ne\emptyset$ and let
$\Gamma$ be the Gelfand transform. Then $\Gamma$ is injective iff $A$ is tractable.
\pf By Lemma A2.4 there is a bijection
$\Omega(A)\to M(A);\,\omega\mapsto{\rm ker}\,\omega$.
We have
$$\eqalign{\ker\Gamma&=\{a\in A\colon\ \ha=0\,\}\cr
&=\{a\in A\colon\ \forall\,\omega\in\OA\ \omega(a)=0\,\}\cr
&=\cap_{\omega\in\OA}\{a\in A\colon\ \omega(a)=0\,\}\cr
&=\cap_{\omega\in\OA}\ker\omega\cr
&=\cap\,M(A),\cr}$$
and the result follows.\eop

\noindent A result concerning tractability, which we state because it is needed in chapter 3, is:

\prop A2.12: Let $A$ be a commutative unital normed algebra and $(\tilde A, i)$
a completion. Then if $\tilde A$ is semisimple, $A$ is tractable.
\pf Suppose $\tilde A$ is semisimple and that  $a\in\bigcap M(A)$. Let
$\omega\in\Omega(\tilde A)$. $i$ is a unital homomorphism so $\omega\circ i
\in\Omega(A)$. So by assumption $(\omega\circ i)(a)=\omega(i(a))=0$. Therefore
$i(a)\in\bigcap M(\tilde A)=(0)$. But $i$ is injective so $a=0$.\eop

\noindent The converse to this is false; see chapter 3 for an example.

\dfn A2.13: Let $A$ be a normed algebra. $K(A):=\cap M(A)$ is the {\it radical} of $A$.

\vskip 10pt
\noindent The last result of this section confirms the result asserted in [{\bf 1}] that factoring
a normed algebra by $\ker\Gamma$ (where it exists) gives a tractable normed
algebra. We prove it in the commutative, unital case; it is applied in the proof of Theorem 2.13.

\def\Gt{\tilde \Gamma}
\def\G{\Gamma}
\prop A2.14: Let $A$ be a commutative unital normed algebra and
$M(A)$ be as above. Then
$A/K(A)$ is a tractable normed algebra.
\pf By the fundamental isomorphism theorem, $\G$ induces an isomorphism $\Gt\map{A/K(A)}\hA$ with norm $\norm\G$. Now $\hA\subseteq C(\OA)$ and we have the evaluation characters
$$\eh\omega\mapto\hA\CN f {f(\omega)}\qquad(\omega\in\OA)$$
It follows that $\hA$ is tractable because if $f\in K(\hA)$ then $\fa\omega\OA \eh\omega(f)=f(\omega)=0$ so $f=0$.

Now $\Gt$ is an algebraic isomorphism so it induces a bijection between the maximal ideals of $A/K(A)$ and $\hA$
given by $J\mapsto \Gt(J)$. $\Gt$ is continuous so for every closed maximal ideal, $J$, of $\hA$, $\Gt^{-1}(J)$ is a closed maximal ideal of $A/K(A)$. Therefore
$$K\left(A/K(A)\right)=\bigcap M\left(A/K(A)\right)\subseteq
\bigcap\left\{\Gt^{-1}(J)\st J\in M(\hA)\right\}=\Gt^{-1}\left(\bigcap M(\hA)\right)=\Gt^{-1}(\{0\})=\{0\}\eqno\square$$

\vfill\eject

\def\a{\alpha}
\def\b{\beta}

\def\BE{{\bf E}}

\font\bigtenrm=cmr10 scaled\magstep {2.5}
\font\bigtenbf=cmbx10 scaled\magstep {1.5}
\noindent {\bigtenbf Appendix 3}
\vskip 7pt
\noindent {\bigtenrm Set Theory}
\vskip 15mm

\noindent The constructions in chapters 5 and 6 require powerful set theory and here we set out the relevant definitions, results,
notations and other conventions. All the results in this chapter are well known and we shall omit all proofs of the more standard ones;
we refer the reader to [{\bf 21}] for these.

\dfn A3.1: Let $(W, \le)$ be a well-ordered set and $\a,\b\in W$ with $\a\le\b$. The {\it segment} $(\a,\b)$ is $\{w\in W\st\a<w<\b\}$ while the {\it interval} $[\a, \b]:=\{w\in W\st \a\le w\le \b\}$. Let $0=\inf W$. An {\it initial segment} of $W$ is a subset of the form $[0, \a)$ or $W$, the former being {\it proper}.
\vskip 10pt
\noindent
The following three lemmas are needed in chapter 5.

\lem A3.2: Let $A$ be an infinite set. Then $\#A=\#(A\times A)$.
\pf A solution to this can be found in [{\bf 21}].\eop

\cor A3.3: Let $A$ be an infinite set and $\GU$ be a well-ordered set with $0<\#\GU\le\#A$. Suppose that $\fa\a\GU \Aa$ is a set with $\#\Aa=\#A$. Then
$$\#\left(\bigcup_{\a\in\GU}\Aa\right)=\#A.$$
\pf Obviously $\#A\le\#\left(\bigcup_{\a\in\GU}\Aa\right)$. Now for $a\in\bigcup_{\a\in\GU}\Aa, \ \te\a\GU\st a\in\Aa$. Also there is a bijection $\phi_\a\map\Aa A$. We have $(\a,\phi_\a(a))\in\GU\times A$ and there is an injection
$\bigcup_{\a\in\GU}\Aa\to\GU\times A$. The result now follows from the Cantor-Schr\"oder-Bernstein theorem since $\#(\GU\times A)\le\#(A\times A)=\#A$.\eop

\cor A3.4: Let $A$ be an infinite ring. Then $\#A=\#A[x]$.
\pf Clearly $\#A\le\#A[x]$. Define
$$\theta\colon A[x]\to\bigcup_{n\in\NN}A^n\,;\;\al=a_0+\cdots+a_nx^n\mapsto
\cases{(a_0,\ldots,a_n)\in A^{n+1} & $a_n\ne0$\cr
	(0)\in A & $\al=0$.\cr}$$
$\theta$ is an injection so $\#A[x]\le\#\left(\bigcup_{n\in\NN}A^n\right)\le\#(\NN\times A)$ as in the proof of Corollary A3.3.
And $\#(\NN\times A)\le\#(A\times A)=\#A$ by Lemma A3.2.\eop

\dfn A3.5: A set, $S$, is {\it transitive} if $\fa s S s\subseteq S$.
\vskip 10pt\noindent
DEFINITION A3.6: The set $\a$ is an {\it ordinal} if $\a$ is transitive and $\a=\emptyset$ or $(\a,\in)$ is well-ordered.

\prop A3.7: A set, $S$, is an ordinal $\iff \exists$ a well-ordering on $S$ such that $\fa s S
[0,s)=s$.

\prop A3.8: Let $\a,\b$ be ordinals. Then $\a\in\b\iff\a\subset\b$.

\dfn A3.9: Whenever $\a,\b$ are ordinals $\a<\b\,:\Leftrightarrow\,\a\subset\b$.

\dfn A3.10: Let $(S, \le), (T, \sqsubseteq)$ be partially ordered sets and $\phi\map S T$ be bijective. $\phi$ is an {\it order isomorphism} if $\fa {s_1,s_2}S\;s_1\le s_2\Leftrightarrow \phi(s_1)\sqsubseteq\phi(s_2)$.

\prop A3.11: If $\phi\map\a\b$ is an order isomorphism of ordinals then $\a=\b$.

\prop A3.12: Let $\BE$ be a non-empty set of ordinals. Then $(\BE, \subseteq)$ is well-ordered.

\prop A3.13: Let $\BE$ be a set of ordinals. Then $\cup\BE$ is an ordinal.

\prop A3.14: Let $\a$ be an ordinal. Then $\a+1:=\a\cup\{\a\}$ is an ordinal.

\dfn A3.15: Let $\BE$ be a set of ordinals and $\b\in\BE$. $\b$ is a {\it non-limit ordinal in $\BE$} if $\te\a\BE\st\b=\a+1$.
Otherwise $\b$ is a {\it limit ordinal in $\BE$}.

\thm A3.16 (The Counting Theorem): Let $W$ be a well-ordered set. Then there is a unique ordinal which is order isomorphic to $W$.

\thm A3.17 (The Principle of Transfinite Induction): Let $W$ be a well-ordered set and $J\subseteq W$ be such that $\fa w W\;[0,w)\subseteq J\,\Rightarrow w\in J$. Then $J=W$.

\thm A3.18 (The Transfinite Recursion Theorem): Let $W$ be a well-ordered set and $X$ a set. Let
$${\cal Q}=\cup_{I\st I\hbox{\sevenrm is an initial segment of }W} {\rm Map}(I, X)$$
Let $f\in\,{\rm Map}({\cal Q}, X)$. Then $\exists\vert\,U\in\,{\rm Map}(W, X)\st\fa a W U(a)=f\left(U|_{[0, a)}\right)$.

\thm A3.19: Let $S\ne\emptyset$. Then there exists a well-ordering on $S$.
\vskip 10pt
\noindent The applications require the existence of an uncountable well-ordered set which has the property that all its proper initial segments are countable. By A3.19.16 there certainly exist uncountable ordinals. By A3.12 there is a unique first uncountable ordinal, which we denote $\omega_1$. It is easy to see that $\omega_1$ has the desired property. The sort of situation in which it becomes useful is illustrated by the result below (which is also used in chapter 6).

\thm A3.20: Let $(X, d)$ be a metric space and $(F_\a)_{0\le\a<\omega_1}$ be a family of closed sets in $X$ such that $\forall\,\a\le\b<\omega_1\;F_\a\subseteq F_\b$. Then $F:=\bigcup_{\a<\omega_1}F_\a$ is closed.
\pf Let $(x_n)_{n\in\NN}\subseteq F$ and $x_n\to x\in X\;(n\to+\infty)$. Now $\fa n \NN\exists\,\a_n<\omega_1\,:\;x_n\in F_{\a_n}$.

Now $\exists\,\a<\omega_1\st\fa n \NN \a_n\le\a$ for if this were not true then $\forall\,\a<\omega_1\; \te n \NN\st \a<\a_n$. Hence $[0, \omega_1)\subseteq\bigcup_{n\in\NN}[0, \a_n)$. But by the definition of $\omega_1$, $[0, \a_n)$ is countable $\fa n \NN$. Therefore $\omega_1=[0, \omega_1)$ is countable since it is contained in a countable union of countable sets. But this is a contradiction.

Therefore $(x_n)\subseteq F_\a$ and since $F_\a$ is closed $x\in F_\a\subseteq F$.\eop

\vfill\eject




\font\bigtenrm=cmr10 scaled\magstep {2.5}
\font\bigtenbf=cmbx10 scaled\magstep {1.5}
\noindent {\bigtenbf References}
\vskip 4mm

\item{[{\bf 1}]} Arens, R. and Hoffman, K. (1956) Algebraic Extension of Normed Algebras. {\sl Proc. Am. Math. Soc.}, 7, 203-210

\vskip 11pt
\item{[{\bf 2}]} Bonsall, F. F. and Duncan, J. (1973) {\sl Complete Normed Algebras.} Germany: Springer-Verlag.

\vskip 11pt
\item{[{\bf 3}]} Bollob\'as, B. (1975) Adjoining Inverses to Commutative Banach Algebras, In: Williamson, J. H. (ed.)
{\sl Algebras in Analysis.} Norwich: Academic Press Inc. (London) Ltd. 256-257

\vskip 11pt
\item{[{\bf 4}]} Lindberg, J. A. (1964) Factorization of Polynomials over Banach Algebras. {\sl Trans. Am. Math. Soc.}, 112, 356-368

\vskip 11pt
\item{[{\bf 5}]} Simmons, G. F. (1963) {\sl Introduction to Topology and Modern Analysis.} Singapore: McGraw-Hill Book Co.

\vskip 11pt
\item{[{\bf 6}]} Jacobson, N. (1996) {\sl Basic Algebra I.} 2nd ed. United States of America: W.H. Freeman and Company.

\vskip 11pt
\item{[{\bf 7}]} Lindberg, J. A. (1964) Algebraic Extensions of Commutative Banach Algebras. {\sl Pacif. J. Math.}, 14, 559-583

\vskip 11pt
\item{[{\bf 8}]} Stout, E. L. (1973) {\sl The Theory of Uniform Algebras.} Tarrytown-on-Hudson, New York: Bogden and Quigley Inc.

\vskip 11pt
\item{[{\bf 9}]} Lindberg, J. A. (1963) On the Completion of Tractable Normed Algebras. {\sl Proc. Am. Math. Soc.}, 14, 319-321

\vskip 11pt
\item{[{\bf 10}]} Cole, B. J. (1968) {\sl One-Point Parts and the Peak-Point Conjecture}, Ph.D. Thesis, Yale University.

\vskip 11pt
\item{[{\bf 11}]} Lindberg, J. A. (1973) Integral Extensions of Commutative Banach Algebras. {\sl Can. J. Math.}, 25, 673-686

\vskip 11pt
\item{[{\bf 12}]} Heuer, G. A. and Lindberg, J. A. (1963) Algebraic Extensions of Continuous Function Algebras. {\sl Proc. Am. Math. Soc.}, 14, 337-342

\vskip 11pt
\item{[{\bf 13}]} Rickart, C. E. (1960) {\sl Banach Algebras.} United States of America: D. Van Nostrand Company, Inc.

\vskip 11pt
\item{[{\bf 14}]} Karahanjan, M. I. (1979) Some Algebraic Characterisations of the Algebra of All Continuous Functions on a
Locally Connected Compactum. {\sl Math. USSR Sb.}, 35, 681-696

\vskip 11pt
\item{[{\bf 15}]} Rotman, J. J. (1979) {\sl An Introduction to Homological Algebra.} London: Academic Press Inc.

\vskip 11pt
\item{[{\bf 16}]} Leibowitz, G. M. (1970) {\sl Lectures on Complex Function Algebras.} United States of America: Scott, Foresman and Company.

\vskip 11pt
\item{[{\bf 17}]} Kelley, J. L. (1955) {\sl General Topology.} New-York: Van Nostrand.

\vskip 11pt
\item{[{\bf 18}]} Feinstein, J. F. (1989) {\sl Derivations from Banach Function Algebras}, Ph.D. Thesis, University of Leeds

\vskip 11pt
\item{[{\bf 19}]} Zame, W. R. (1984) Covering Spaces and the Galois Theory of Commutative Banach Algebras. {\sl J. Funct. Anal.}, 27, 151-171

\vskip 11pt
\item{[{\bf 20}]} Murphy, G. J. (1990) {\sl C*-Algebras and Operator Theory.} Boston: Academic Press.

\vskip 11pt
\item{[{\bf 21}]} Halmos, P. R. (1970) {\sl Naive Set Theory.} New-York: Springer-Verlag.

\eject
\end